\documentclass[letterpaper]{IEEEtran}  

\usepackage{times} 
\usepackage{amsmath} 
\usepackage{amssymb}
\usepackage{amsfonts}  
\usepackage{latexsym}
\usepackage{theorem}
\usepackage{graphicx} 
\usepackage{dsfont}
\usepackage{mathtools}
\usepackage{subfigure}
\usepackage{lscape}
\usepackage{color}
\usepackage{float}
\usepackage{enumerate}
\usepackage[percent]{overpic}

	\newtheorem{theorem}{Theorem}

	\newtheorem{assumption}{Assumption}
	
	\newtheorem{corollary}[theorem]{Corollary}
	\newtheorem{definition}{Definition}
	{\theorembodyfont{\rmfamily} }
	\newtheorem{lemma}[theorem]{Lemma}

	\newtheorem{proposition}[theorem]{Proposition}

	\allowdisplaybreaks
	\IEEEoverridecommandlockouts
	
	\DeclareMathOperator*{\argmin}{arg\,min}
	\DeclareMathOperator*{\argmax}{arg\,max}
	\title{Fast Convergence Rates for Distributed Non-Bayesian Learning}

	\author{Angelia Nedi\'{c}, Alex Olshevsky and C\'{e}sar A. Uribe$^*$	\thanks{$^*$Corresponding author.}
		\thanks{A. Nedi\'{c} (\textit{angelia.nedich@asu.edu}) is with the ECEE Department, Arizona State University. A. Olshevsky (\textit{alexols@bu.edu}) is with the Department of ECE and Division of Systems Engineering, Boston University. C.A. Uribe (\textit{cauribe2@illinois.edu}) is with the Coordinated Science Laboratory, University of Illinois.
			This research is supported partially by the National Science Foundation under
			grants no.\ CCF-1017564 and no.\ CMMI-1463262 and by the Office of Naval Research under grant
			no.\ N00014-12-1-0998.}  
	}
	
\begin{document}
		\maketitle
		\begin{abstract}
			We consider the problem of \textit{distributed learning}, where a network of agents collectively aim to agree on a hypothesis that best explains a set of distributed observations of conditionally independent random processes. We propose a distributed algorithm and establish consistency, as well as a non-asymptotic, explicit and geometric convergence rate for the concentration of the beliefs around the set of optimal hypotheses. Additionally, if the agents interact over static networks, we provide an improved learning protocol with better scalability with respect to the number of nodes in the network.  
		\end{abstract}
		
		
		\section{Introduction}
		
Large numbers of interconnected components add to the complexity of engineering systems. Developing models and tools for the analysis of such distributed systems is necessary, not only from the engineering point of view but for effective decision-making and policy design. For example, the control of autonomous vehicle for exploration, rescue, and surveillance depends on the coordination abilities of fleets of robots; each robot should make decisions based on local information and limited communications. Power networks (e.g. the electric grid) need several generating and consuming stations to coordinate offer and demand to improve efficiency. In traffic control, the goal is to distributively avoid jams and to increase traffic flow based on limited infrastructure (e.g. roads). Economic systems need modeling, estimation and control of markets at the micro and macroeconomic scales. Market dynamics depend on several agents influencing the system, each of which might have conflicting goals. In telecommunication networks, several stations need to communicate over non-perfect channels to optimize information transmission. The control of industrial processes requires communication and coordination between different parts of the process in hazardous environments. The modeling and control of ecological systems requires the analysis of several actors interacting with each other, subject to changing environments.  

Traditional approaches for the design of distributed inference algorithms, for inherently distributed systems, assume a fusion center exists. The fusion center gathers all the information and makes centralized decisions \cite{gub93,zhu05,vis97,sun04}. Nonetheless, communication constraints, limited memory and lack of physical accessibility to certain measurements hinders this task. Therefore, it is necessary to develop algorithmic protocols that take into account such constraints and use only locally available information. Although many results on these themes have appeared in recent years, the study of distributed decision-making and computation traces back to classic papers from the 70s and 80s~\cite{aum76,bor82,tsi84,gen86,coo90,deg74,gil93}.

In~\cite{ace11}, the authors describe results on learning in social networks based on computing posterior distributions using Bayes' rule. That is, given some assumed prior knowledge and new observations, an agent computes a posterior based on likelihood models, see~\cite{mue13}. Nevertheless, a fully Bayesian approach might not be possible because full knowledge of the network structure, or other agents' likelihood models, need not be available~\cite{gal03,mos10}. Other authors showed that non-Bayesian methods can be used in learning task as well~\cite{eps10,jad12,jad13,sha13}. In this case, agents are assumed to be boundedly rational (i.e. fail to aggregate information in a fully Bayesian manner \cite{gol10}). They repeatedly communicate with others and use naive approaches to aggregate information.

 Several groundbreaking papers have described distributed methods to achieve global behaviors by repeatedly aggregating local information without complete knowledge of the network~\cite{jad12,rah10,ala04,olf06}. For example, in distributed hypothesis testing using belief propagation, convergence and its dependence on the communication structure were shown~\cite{ala04}. Later, extensions to finite capacity channels, packet losses, delayed communications and tracking were developed~\cite{sal06,rah07}. In~\cite{rah10}, the authors proved convergence in probability, the asymptotic normality of the distributed estimation and provided conditions under which the distributed estimation is as good as a centralized one. Later in~\cite{jad12}, the almost sure convergence of a non-Bayesian rule based on arithmetic mean was shown for fixed topology graphs. Extensions to information heterogeneity and asymptotic convergence rates have been derived as well~\cite{jad13}. Following~\cite{jad12}, other methods to aggregate Bayes estimates in a network have been explored. In \cite{ban14}, geometric means are used for fixed topologies as well, however, the consensus and learning steps are separated. The work in~\cite{qip15} extends the results of~\cite{jad12} to time-varying undirected graphs. In~\cite{sha13}, local exponential rates of convergence for undirected gossip-like graphs are studied. The authors in~\cite{lal14b,sha14,qip11,qip15} proposed a non-Bayesian learning algorithm where a local Bayes' update is followed by a consensus step. In~\cite{lal14b}, convergence result for fixed graphs is provided and large deviation convergence rates are given, proving the existence of a random time after which the beliefs will concentrate exponentially fast. In~\cite{sha14}, similar probabilistic bounds for the rate of convergence are derived for fixed graphs and comparisons with the centralized version of the learning rule are provided. Other variations of the non-Bayesian approach have been proposed for continuum set of hypotheses \cite{ned16d}, weakly connected graphs \cite{sal16}, bisection search algorithms \cite{tsi16}, transmission node failures ~\cite{su16,su16b,su16c} and time-varying graphs \cite{ned14,ned15b,ned16}. See \cite{bar13,ned16c} for an extended literature review.

In this paper, we consider a network of agents, where each agent repeatedly receives information from its neighbors as well as private signals from an external source. The private signals are realizations of a random variable with an \textit{unknown} distribution. The agents would like to collectively agree on a hypothesis (distribution) that best explains the data observed by all nodes/agents. We focus on the case where agents might have inconsistent hypotheses, in the sense that, the hypotheses that best describe private observations need not be the same as the hypotheses that best describe the aggregated set of observations of all agents.     

The contributions of this paper are: first, we propose and motivate a novel distributed non-Bayesian learning rule. We derive the proposed algorithm as the solution of a natural extension of the variational representation of Bayes' updates in a distributed setting. This characterizes a general family of distributed non-Bayesian learning protocols. We show that existing protocols are instances of this general family algorithms. Additionally, we show that the proposed protocol allows the network to learn the set of hypotheses that best explain the data collected by all the nodes (i.e. consistency). We also provide a geometric, non-asymptotic, and explicit characterization of the convergence rate, which immediately leads to finite-time bounds that scale intelligibly with the number of nodes for general time-varying undirected graphs. Finally, we propose and analyze a new protocol for arbitrary fixed undirected graphs that scales better than previous algorithms with respect to the number of agents in the network.

Simultaneous and independent works obtained results which overlap with ours \cite{sha14,lal14b}. Specifically, in \cite{sha14} the authors proposed a variant of a distributed learning algorithm, a similar convergence rate was obtained. Consistency and asymptotic rates were provided in \cite{lal14b} for another class of non-Bayesian learning. Moreover, specific instances of the problem studied in this work have been considered in the context of distributed parameter estimation \cite{sal06, rah10}. We note that, relative to these simultaneous papers, our results are more general in the sense that they allow time-varying networks and allow nodes to have conflicting hypotheses, none of which matches the distribution of the observations. Furthermore, in the case of fixed undirected graphs, we propose an update rule which involves an additional register of memory in each node to obtain a more graceful scaling with the number of nodes in the network. Section \ref{general} provides a more detailed comparison with the mentioned papers.

This paper is organized as follows. 
In Section~\ref{model}, we describe the problem and main results. In Section~\ref{general}, we introduce a general class of distributed non-Bayesian learning rules and provide comparisons with recent literature. 
In Section~\ref{consistency}, we analyze the consistency of the information aggregation and estimation models, while in Section~\ref{rates} we prove a non-asymptotic convergence rate for the concentration of the beliefs generated by the proposed algorithm for time-varying graphs. In Section~\ref{linear_rates}, we show the convergence time improvement for a new protocol for fixed undirected graphs. Section~\ref{special} develops the application of the proposed methods for the problem of distributed source localization.  Conclusions and future work directions are discussed in Section~\ref{conclusions}.
		
		\subsubsection*{Notation}
		We use upper case letters to represent random variables (e.g. $X_k$), 
		and the corresponding lower case letters for their realizations (e.g. $x_k$). 
		We write $[A]_{ij}$ or $A_{ij}$ to denote the entry of the matrix $A$ in the $i$-th row and $j$-th column.
		We write $A'$ for the transpose of a matrix $A$ and $x'$ for the transpose of a vector $x$. 
		We use $I_n$ for the identity matrix of size $n$ by $n$. 
		Bold letters represent vectors which are assumed to be column vectors unless specified otherwise. 
		The $i$-th entry of a vector will be denoted by a superscript $i$, i.e.,  
		$\boldsymbol{x}_k = [x_k^1,\hdots,x_k^n]'$.
		We write $\boldsymbol{1}_n$ to denote the all-ones vector of size $n$. 
		For a sequence of matrices $\{A_k\}$, we let  
		$A_{k_f:k_i} \triangleq A_{k_f} \cdots A_{k_i+1} A_{k_i}$ for all $k_f \geq k_i \geq 0$. 
		We abbreviate terminology \textit{almost surely }by {\it a.s.} and \textit{independent identically distributed} by {\it i.i.d.}. 
		In general, when referring an agent $i$ we will use superscripts and when referring to a time instant $k$ we will use subscripts.

		\section{Problem Setup and Main Results}\label{model}
		
		Consider a group of $n$ agents, indexed by $1,2,\ldots,n$, 
		each having observations of conditionally independent random processes, at discrete time steps $k=1, 2, 3, \ldots$. 
		Specifically, agent $i$ observes the random variables $S_1^i, S_2^i, \ldots,$ which are i.i.d. and distributed according  to an unknown
		probability distribution $f^i$. The output space of the random variables $S_k^i$ is a finite set which we will denote by
		${\mathcal{S}}^i$. For convenience, we stack up all the $S_k^i$ into a vector denoted as $\boldsymbol{S}_k$. 
		Then, $\boldsymbol{S}_k$ is an i.i.d.\ vector taking values in 
		${\boldsymbol{\mathcal{S}} = \prod_{i=1}^n \mathcal{S}^i}$ and distributed as ${\boldsymbol{f}=\prod_{i=1}^n f^i}$. Furthermore, each agent $i$ has a family of probability distributions $\{\ell^i(\cdot | \theta)\}$ parametrized by a finite set
		$\Theta = \{ \theta_1, \theta_2, \ldots, \theta_m\}$ with $m$ elements. 
		One can think of $\Theta$ as a set of hypotheses and $\ell^i( \cdot | \theta)$ as the probability distribution that would be seen by agent $i$ if hypothesis $\theta$ were true. We do not require that there exists
		$\theta \in \Theta$ with $\ell^i(\cdot | \theta) = f^i$ almost everywhere for all $i=1, \ldots, n$; in other words, there may not be a hypothesis that matches the observations made by the nodes. Rather, the objective of all agents is to 
		agree on a subset of $\Theta$ that best fits all the observations in the network. 
		
\begin{figure}[ht]
			\centering
			\begin{overpic}[width=0.3\textwidth]{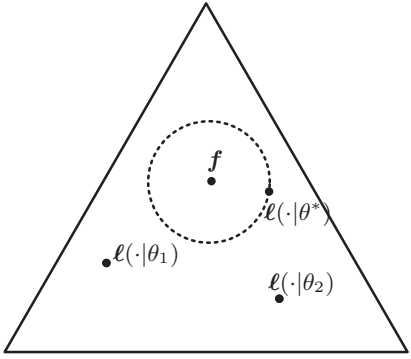}
				\put(50,46){{\small $\boldsymbol{f}$}}
				\put(64,33){{\small $\boldsymbol{\ell}(\cdot|\theta^*)$}}
				\put(27,23){{\small $\boldsymbol{\ell}(\cdot|\theta_1)$}}
				\put(65,16){{\small $\boldsymbol{\ell}(\cdot|\theta_2)$}}
			\end{overpic} 
			\caption{Geometric interpretation of the learning objective. The triangle represents the simplex of all possible probability distributions of $\boldsymbol{S}_k$. The point
				$\boldsymbol{f}$ is the actual distribution of $\boldsymbol{S}_k$. The goal of the network of agents is to learn the hypothesis $\theta^*$ that best describes its observations, which corresponds to the distribution $\boldsymbol{\ell}(\cdot|\theta^*)$ (the closest to the distribution $\boldsymbol{f}$).}
			\label{triangle}
		\end{figure}

		Formally, this setup describes the scenario where the group of agents collectively tries to solve the following optimization problem
		\begin{align}\label{opt_problem}
		\min_{\theta \in \Theta} F(\theta) & \triangleq D_{KL}\left(\boldsymbol{f}\|\boldsymbol{\ell}\left(\cdot | \theta\right)\right) \ \ \  \nonumber \\
		& = \sum\limits_{i=1}^n D_{KL}\left(f^i\|\ell^i\left(\cdot|\theta\right)\right) 
		\end{align}
		where $D_{KL}\left(f^i\|\ell^i\left(\cdot|\theta\right)\right)$ is Kullback-Leibler (KL) divergence between the distribution of $S_k^i$ and $\ell^i( \cdot | \theta)$. The distributions $f^i$'s are unknown, therefore the agents try to ``learn" the solution to this optimization problem based on local observations and interactions, see Figure~\ref{triangle}. 
		
		Consider for example a group of two agents, labeled by 1 and 2, such that $S_k^i \sim \mathcal{N}(i,1)$, which is equivalent to $S_k^i = i + W_k^i$ where $W_k^i \sim \mathcal{N}(0,1)$ is a zero mean Gaussian process with unitary standard deviation. They want to correctly identify the parameter $\theta^*$ out of three possible hypotheses $\Theta = \{\theta_1,\theta_2,\theta_3\}$ where the likelihood models of the agents are:
		$\ell^1(s^1|\theta_1) =  \phi(s^1- 0.5)$, $\ell^2(s^2|\theta_1) =  \phi(s^2)$, $\ell^1(s^1|\theta_2) = \phi(s^1- 1.5)$,  $\ell^2(s^2|\theta_2) =  \phi(s^2- 2.5)$, $\ell^1(s^1|\theta_3) =  \phi(s^1)$, $\ell^1(s^2|\theta_3) =  \phi(s^2 - 1.5)$ where  $\phi(x) = \exp(-\frac{1}{2}x^2) / \sqrt{2\pi}$ is the probability density function of the standard normal distribution.
		In this scenario, agent 1 alone would not be able to differentiate between $\theta_1$ and $\theta_2$ and agent 2 cannot differentiate between $\theta_2$ and $\theta_3$ given that they are at the same distance to the true distribution of the observations. Nonetheless, when they interact with each other the solution to the proposed optimization problem is $\theta^* = \theta_2$.
		
		\subsection{Proposed Learning Algorithms}
		
		Probability distributions over the hypothesis set $\Theta$ will be referred as {\em beliefs}. Every
		agent $i$ has an initial belief $\mu_0^i$, which we often refer to as its \textit{prior distribution} or \textit{prior belief}. We will be studying the dynamics wherein agents exchange beliefs with their neighbors over some communication network, with the effect 
		that over time these beliefs concentrate on the ``best'' choice of hypotheses.
		Each agent $i$ generates a new belief for time $k+1$, which we will denote by $\mu_{k+1}^i$, based on its 
		current belief $\mu_k^i$, an observation $s_{k+1}^{i}$ of the random variable $S^i_{k+1}$, and the current beliefs of its neighbors $\mu_k^j$ with $j \neq i$.  We propose two algorithms for the generation of the new belief $\mu_{k+1}^i$: a generic rule for undirected time-varying graphs and a special rule for static graphs. We show that the proposed update rules generate a sequence of beliefs that sequentially approaches a solution to the optimization problem in \eqref{opt_problem}. 
		
		We consider the following rule for general undirected time-varying graphs: for each $\theta \in \Theta$,
		\begin{align}\label{non_bayes_distributed}
		\mu_{k+1}^i(\theta) & = \frac{1}{Z_{k+1}^i}\prod_{j=1}^n\mu_k^j(\theta)^{[A_k]_{ij}}\ell^i(s_{k+1}^i|\theta)^{\beta^i_k}
		\end{align} 
		where $Z_{k+1}^{i}$ is a normalization factor to make the beliefs a probability distribution, i.e.,
		\begin{align*}
		Z_{k+1}^i & = \sum_{p=1}^m \prod_{j=1}^n \mu_k^j\left(\theta_p\right)^{[A_k]_{ij}}\ell^i(s_{k+1}^i|\theta_p)^{\beta^i_k}
		\end{align*}
		where the $A_k$ is a non-negative matrix of ``weights", which is compliant with the connectivity structure of the underlying communication network. The network at each time instant $k$ is modeled as a graph $\mathcal{G}_k$ is composed 
		by a node set $V = \left\{1,2, \hdots, n \right\}$ and a set $E_k$ of undirected links. The variable $\beta^i_k$ is a stationary Bernoulli random process with mean $q^i$, which indicates if an agent obtained a new realization of $S^i_{k+1}$. 
		Specifically, $\beta^i_k=1$ indicates that agent $i$ obtained a new observation, while $\beta^i_k=0$ indicates that it did not.
		
		For static undirected graphs, we propose a new belief update rule with one-step memory as follows: for each $\theta$ in $\Theta$
		\begin{align}\label{linear_bayes}
		\mu_{k+1}^i(\theta) & 
		= \frac{1}{\widetilde{Z}_{k+1}^i} 
		\frac{ \prod\limits_{j=1}^{n} \mu_k^j(\theta)^{\left(1+\sigma\right)
				\bar{A}_{ij} } \ell^i (s^i_{k+1}|\theta)^{\beta^i_k}  }{\prod\limits_{j=1}^{n} \left( \mu_{k-1}^j(\theta)\ell^j (s^{j}_{k}|\theta)^{\beta^j_{k-1}}\right) ^{\sigma\bar{A}_{ij} }}
		\end{align}
		where $\widetilde{Z}^i_{k+1}$ is the corresponding normalization factor given by
		\begin{align*}
		\widetilde{Z}_{k+1}^i & 
		= \sum\limits_{p=1}^{m}\frac{ \prod\limits_{j=1}^{n} \mu_k^j\left(\theta_p\right)^{\left(1+\sigma\right)
				\bar{A}_{ij} } \ell^i \left(s^i_{k+1}|\theta_p\right) ^{\beta^i_k} }{\prod\limits_{j=1}^{n} \left( \mu_{k-1}^j(\theta_p) \ell^j (s^{j}_{k}|\theta_p)^{\beta^j_{k-1}}\right) ^{\sigma\bar{A}_{ij} }}
		\end{align*}
		where $\bar{A}$ is a specifically chosen matrix (called the \textit{lazy Metropolis matrix}) and $\sigma$ a constant to be set later. We initialize $\mu_{-1}^i (\theta)$ to be equal to $\mu_0^i (\theta)$ for all $i=1, \ldots,n$ and $\theta \in \Theta$.  We will show that this update rule generates a sequence of beliefs that concentrate at a rate 
		a factor of $n$ faster than the previous results.  Note that the update rule described in Eq.~\eqref{linear_bayes} requires the communication of the product of the beliefs and likelihood functions and an additional memory since the beliefs at time $k+1$ depends on the beliefs a time $k$ and at time $k-1$.
		
		Section \ref{general} will motivate the choice of the update rules. They can be interpreted as natural generalizations of the variational representation of the Bayes update rule for the distributed learning setting.
		
		\subsection{Assumptions and Definitions}
		
		We will list a sequence of assumptions about the underlying communication graph and the family of parametrized likelihood models. They will guarantee the desired convergence properties.  
		
		For the first class of update rules described in Eq.~\eqref{non_bayes_distributed}, 
		we assume the following structure for the sequence of communication graphs $\{\mathcal{G}_k\}$.
		\begin{assumption}\label{assum:graph}
			The graph sequence $\{ \mathcal{G}_k \}$ and the matrix sequence $\{A_k\}$ are such that:
			\begin{enumerate}[(a)]
				\item $A_k$ is doubly-stochastic with $[A_k]_{ij} > 0$ if $( i,j)\in E_k$.
				\item  If $(i,j) \notin E_k$ for some $i \neq j$ then $A_{ij}=0$.
				\item $A_k$ has positive diagonal entries, $[A_k]_{ii}>0$ for all $i=1,\ldots,n.$
				\item If $[A_k]_{ij}>0$, then $[A_k]_{ij} \ge \eta$ for some positive constant $\eta$.
				\item $\left\{\mathcal{G}_k\right\}$ is $B$-strongly connected, i.e.,
				there is an integer $B\ge 1$ such that the graph 
				$\left\{V,\bigcup_{i=kB}^{\left(k+1\right)B-1}E_i\right\}$ is strongly connected for all $k \geq 0$.
			\end{enumerate}
		\end{assumption} 
		
		Assumption~\ref{assum:graph}(a) and Assumption \ref{assum:graph}(b) characterize the communication between agents. If two agents can exchange information at a certain time instant $k$, the underlying communication graph will have an edge between the corresponding nodes. This also implies a positive weighting of the information shared. 
		The graph sequence $\{ \mathcal{G}_k \}$ and the matrix sequence $\{A_k\}$ define a corresponding inhomogeneous Markov Chain with transition probabilities $A_k$. 
		Assumption~\ref{assum:graph}(c) guarantees the aperiodicity of this Markov Chain. Additionally, Assumptions~\ref{assum:graph}(d) and~\ref{assum:graph}(e) guarantee that this Markov chain is ergodic by ensuring there is sufficient connectivity and that the entries of $A_k$ do not vanish. Assumption~\ref{assum:graph} is common in distributed optimization and consensus literature \cite{ned09,ned13}. It guarantees convergence of the associated Markov Chain and defines bounds on relevant eigenvalues in terms of the number of agents.
		
		There are several ways to construct a set of weights satisfying Assumption~\ref{assum:graph}. For example, one can consider a lazy Metropolis (stochastic) matrix of the form 
		$\bar{A}_k = \frac{1}{2}I_n + \frac{1}{2}\hat{A}_k$, where $I_n$ is the identity matrix and $
		\hat{A}_k$ is a stochastic matrix whose off-diagonal entries satisfy
		\begin{align*}
		[ \hat{A}_k] _{ij} = \left\{  \begin{array}{l l}
		\frac{1}{\max \left\{ d^i_k+1,d^j_k+1 \right\}  }, & \quad \text{if $(i,j) \in E_k$ }\\
		0 ,& \quad \text{if $(i,j) \notin E_k$}
		\end{array} \right. 
		\end{align*}
		where $d^i_k$ is the degree (the number of neighbors) of node $i$ at time $k$. Note that the lazy Metropolis weights require undirected communications since each weight $[\hat{A}_k]_{ij}$ depends on the degree of both agent $i$ and agent $j$. Thus, we will require that agents share their beliefs as well as their degree, which means exchanging $m+1$ numbers at each time step.
		
	Analogous to Assumption~\ref{assum:graph}, we use the following assumption when the interaction between the agents happens over static graphs with the update rule described in Eq.~\eqref{linear_bayes}.
		\begin{assumption}\label{assum_linear}
			The graph sequence $\{ {\cal G}_k\}$ is static (i.e. $\mathcal{G}_k = \mathcal{G}$ for all $k$) and undirected and the weight matrix
			$\bar{A}$ is a lazy Metropolis matrix, defined by
			\begin{align*}
			\bar{A} = \frac{1}{2}I_n + \frac{1}{2}\hat{A}
			\end{align*}
			where $\hat{A}$ is the \textit{Metropolis matrix}, which is the unique stochastic matrix whose off-diagonal entries satisfy
			\begin{align*}
			\hat{A}_{ij} = \left\{  \begin{array}{l l}
			\frac{1}{\max \left\{ d^i+1,d^j+1 \right\} }, & \quad \text{if $(i,j) \in E$ }\\
			0, & \quad \text{if $(i,j) \notin E$}
			\end{array} \right. 
			\end{align*}
			with $d^i$ being the degree of the node $i$ (i.e., the number of neighbors of $i$ in the graph).
		\end{assumption}
		
		Next, we provide three important definitions that we use in the sequel to describe some learning-related quantities. 
		\begin{definition}\label{confi}
			The \textit{group confidence} of a nonempty subset ${{W \subseteq V}}$ of agents is given by
			\begin{align*}
			\mathsf{C}^W_{\boldsymbol{q}}(\theta)=-\sum\limits_{i \in W} q^iD_{KL}\left(f^i\|\ell^i\left(\cdot|\theta\right)\right)\qquad\hbox{for all }\theta\in\Theta,
			\end{align*}
			where $q^i$ is the mean-value of the i.i.d.\ Bernoulli variable $\beta_k^i$ characterizing the availability of measurements for agent $i$. If $W=V$, we simply write $\mathsf{C}_{\boldsymbol{q}}$.
		\end{definition}
		
		The group confidence provides a way to quantify the quality of a hypothesis from the perspective of a subset of the agents. 
		The quality of a hypothesis for individual agents is weighted by the mean  
		of the i.i.d. Bernoulli process governing the availability of observations. 
		
		\begin{definition}
			Two distinct hypotheses $\theta_i$ and $\theta_j$ 
			are said to be $W$-observationally equivalent if 
			$\mathsf{C}^W_{\boldsymbol{q}}\left(\theta_i\right) = \mathsf{C}^W_{\boldsymbol{q}}\left(\theta_j\right)$.
		\end{definition}
		
This definition extends the idea of observational equivalence introduced in \cite{jad12}. Group observational equivalence provides a general definition where a group of agents can not differentiate between two hypotheses even if their corresponding likelihood models are not the same.
		
		Finally, we introduce the optimal set of hypotheses as the set with the maximum group confidence. 
		\begin{definition}\label{def:opt}
			The optimal hypothesis set is defined as  
			${\Theta^* = \argmax\limits_{\theta \in \Theta} \mathsf{C}_{\boldsymbol{q}}(\theta)}$, 
			and the confidence of the optimal hypothesis set is denoted as $\mathsf{C}^*_{\boldsymbol{q}}$, i.e.,
			$\mathsf{C}^*_{\boldsymbol{q}} = \mathsf{C}_{\boldsymbol{q}}(\theta^*)$ for  $\theta^*\in\Theta^*$.
		\end{definition}
		
The optimal set is always nonempty, and we assume it is a strict subset of $\Theta$ to avoid the trivial case where all hypotheses are observationally equivalent. This holds if there is a unique true state, $\hat{\theta} \in \Theta$, such that each agent $i$ sees distributions generated according to $f^i = \ell^i(\cdot | \hat{\theta})$, and $\Theta$ contains other hypotheses besides $\hat{\theta}$. 

Informally, we will refer to our assumptions above as describing a setup with \textit{conflicting models}; by this, we mean that the hypothesis which best describes the observations of agent $i$ (i.e., the hypothesis $\theta$ which minimizes $D_{\rm KL}(f^i\| \ell^i(\cdot|\theta))$) may not be the hypothesis which best describes the observations of a different agent, and may in fact not belong to the optimal set $\Theta^*$.

		We will further require the following assumption on the agents' prior distributions and likelihood functions. 
		The first of these is sometimes referred to as the Zero Probability Property~\cite{gen86}. 
		
		\begin{assumption}\label{assum:init}
			For all agents $i=1,\ldots,n$,
			\begin{enumerate}[(a)]
				\item
				The set $\hat{\Theta}^* = \cap_{i=1}^n\Theta^{*i}$ is nonempty, where 
				$\Theta^{*i}\subseteq \Theta^*$ is the subset of optimal hypotheses with positive initial beliefs for agent $i$, i.e., 
				$\mu_{0}^i(\theta)>0$ for all $\theta \in \Theta^{*i}$ and $\mu_{0}^i(\theta)=0$ for all 
				$\theta \in \Theta^*\setminus\Theta^{*i}$.
				\item 
				The support of the true distribution of the observations is contained in the support of the likelihood models for all hypothesis, i.e., there exists an $\alpha >0$ such that if $f^i\left(s^i\right)>0$ then  
				$\ell^i\left(s^i | \theta \right)  > \alpha$ for all $\theta \in \Theta$.
			\end{enumerate}
		\end{assumption}
		
Uniform prior beliefs satisfy the Assumption~\ref{assum:init}(a), which is a reasonable assumption if there is no initial information about the hypotheses quality. In Eq.~\eqref{non_bayes_distributed}, if $\mu^i_k(\theta)=0$ for some hypothesis $\theta$ and for some agent $i$, at some instance $k$, then all beliefs of all agents will eventually become zero at that hypothesis. Assumption~\ref{assum:init}(a) removes the undesired effects of this property which could lead to the inability to learn. In addition, Assumption~\ref{assum:init}(b) guarantees the sub-Gaussian behavior of the observed random variables. Specifically, the derived convergence rates use results from the measure concentration of random variables. In the most common setting, the random variables must have a sub-Gaussian or sub-exponential behavior \cite{bou13}. 

 		\subsection{Results}
		
		We now state our first result; we show that the dynamics in Eq.~\eqref{non_bayes_distributed}
		concentrates the beliefs on the optimal set $\Theta^*$,
		which is precisely the set that best describes the observations. This theorem will be proven in Section~\ref{consistency}.
		
		\begin{theorem}\label{proof_consistency} Under
			Assumptions~\ref{assum:graph} and~\ref{assum:init}, the update rule of Eq.~(\ref{non_bayes_distributed}) 
			has the following property:
			\begin{align*}
			\lim_{k \to \infty} \mu_k^i(\theta) =  0 \quad a.s. 
			\quad\hbox{for all }\theta \notin \hat\Theta^*, ~i = 1, \ldots, n.
			\end{align*} 
		\end{theorem}
		
		Our results regarding the non-asymptotic explicit convergence rate of the update rules in 
		Eq.~(\ref{non_bayes_distributed}) and Eq.~\eqref{linear_bayes} are given in Theorem~\ref{teo2} and Theorem~\ref{thm_lineal}, while their proofs are provided in Section~\ref{rates} and Section~\ref{linear_rates}, respectively.
		
		\begin{theorem}\label{teo2} 
			Let Assumptions~\ref{assum:graph} and~\ref{assum:init} hold and let $\rho\in(0,1)$. The update rule of Eq.~(\ref{non_bayes_distributed}) has the following property:
			there is an integer $\boldsymbol{N}(\rho)$ such that, with probability $1 -\rho$, 
			for all $k\ge \boldsymbol{N}(\rho)$ and for all $\theta_v\notin\Theta^*$, we have
			$$\mu_{k}^i(\theta_v) \leq \exp\left( -\frac{k}{2}\gamma_2+ \gamma_1^i\right)
			\quad\hbox{for all } i = 1, \ldots, n$$ 
			where
			\begin{align*}
			\boldsymbol{N}(\rho)&
			\triangleq \left\lceil\frac{1 }{\gamma_2^2} 8 \left( \log\alpha\right) ^2  \log\frac{1}{\rho} \right\rceil\\
			\gamma_1^i & \triangleq \max_{\substack{\theta_w\in \hat{\Theta}^* \\\theta_v \notin \Theta^*} } \left\{
			\max_i \log \frac{\mu_0^i(\theta_v)}{\mu_0^i(\theta_w)}
			 \right\} + \frac{12 \log n}{1 - \lambda} \log\frac{1}{\alpha} \\
			\gamma_2 &\triangleq \frac{1}{n}\, \min_ {\theta_v\notin\Theta^*}
			\left(\mathsf{C}^*_{\boldsymbol{q}} -\mathsf{C}_{\boldsymbol{q}}(\theta_v)\right)
			\end{align*}
			with $\alpha$ from Assumption~\ref{assum:init}(b), $\eta$ from Assumption~\ref{assum:graph}(d)
			and $\lambda$ given by:
			\begin{align*}
			\lambda = \left(1-\frac{\eta}{4n^2} \right)^{\frac{1}{B}}.
			\end{align*} 
			If each $A_k$ is the 
			lazy Metropolis matrix associated with $\mathcal{G}_k$ and $B=1$, then
			\begin{align*}
			\lambda = 1- \frac{1}{\mathcal{O}(n^2) }.
			\end{align*}
		\end{theorem} 
		
		In words, the belief of each agent on any hypothesis outside the optimal set decays at a network-independent rate which scales with the constant $\gamma_2$, which is the
		average Kullback-Leibler divergence to the next best hypothesis. However, there is a transient due to the $\gamma_1^i$ term (since the bound of 
		Theorem \ref{teo2} is not below $1$ until $k \geq 2 \gamma_1^i/\gamma_2$), and the size of this transient depends on the network and the
		number of nodes through the constant $\lambda$.
		
		Observe that the term $\gamma_1^i$ represents the influence of the initial beliefs as well as the mixing properties of the graph. 
		If all agents use uniform initial beliefs, i.e., $\mu_0^i\equiv 1 / |\Theta|$, then the effect of the initial beliefs is zero and 
		$\gamma_1^i$ reduces to 
		\[\gamma_1^i =
		 \frac{12 \log n}{1 - \lambda} \log\frac{1}{\alpha}\]
		where the constant $\lambda$ may be thought of as the ``time to ergodicity'' of the inhomogeneous Markov Chain associated with the matrix sequence $A_k$. On the other hand, if one can start with an informative prior where $\mu_0^i(\theta^*) > \mu_0^i(\theta)$, the influence of the initial beliefs will be a negative term, effectively reducing the transient time.
		
		Our next result shows the belief concentration rate for the update rule described in Eq.~\eqref{linear_bayes}.
		
		\begin{theorem}\label{thm_lineal} 
			Let Assumptions~\ref{assum_linear} and~\ref{assum:init} hold and let $\rho\in(0,1)$. Furthermore let $U \geq n$ and let $\sigma = 1 - 2/(9U+1)$.	Then, the update rule of Eq.~\eqref{linear_bayes} with this $\sigma$, uniform initial beliefs with the condition $\mu^i_{-1}(\theta)=\mu^i_{0}(\theta)$ and $\beta_{-1}^i$ fixed to zero, has the following property:
			there is an integer $\boldsymbol{N}(\rho)$ such that,
			with probability $1 -\rho$, 
			for all $k\ge \boldsymbol{N}(\rho)$ and for all $\theta_v\notin\Theta^*$,
			it holds that
			$$\mu_{k}^i(\theta_v) \leq \exp\left( -\frac{k}{2}\gamma_2+ \gamma_1^i\right)
			\quad\hbox{for all  }i = 1, \ldots, n,$$ 
			where
			\begin{align*}
			\boldsymbol{N}(\rho) &
			\triangleq \left\lceil\frac{1}{\gamma_2^2}48 \left( \log\alpha\right) ^2\log\left(\frac{1}{\rho} \right)  \right\rceil\\
			\gamma_1^i & \triangleq 
			\frac{4 \log n}{1 - \lambda} \log\frac{1}{\alpha}  \\
			\gamma_2 &\triangleq \frac{1}{n}\, \min_ {\theta_v\notin\Theta^*}
			\left(\mathsf{C}^*_{\boldsymbol{q}} -\mathsf{C}_{\boldsymbol{q}}(\theta_v)\right)
			\end{align*}
			with $\alpha$ from Assumption~\ref{assum:init}(b) and $\lambda = 1 - \frac{1}{18U}$.
		\end{theorem} 
		
 Note that the beliefs for $k=-1$ and $k=0$ are defined equal. Additionally, we assume there is no observation available for time $0$, this holds if we assume $\beta_{-1}^i=0$ with any realization of $S_{0}^i$. 
 
 The bound of Theorem 3 is an improvement by a factor of $n$ compared to the bounds of Theorem \ref{teo2}. In a network of $n$ agents where $\alpha$, $\rho$ and $\gamma_2$ are treated like constants with respect to the number of agents, we require at least $\mathcal{O}(n\log n)$ iterations for the beliefs on the incorrect hypotheses to be below certain small value epsilon (assuming $U$ is within a constant factor of $n$). Following the results of \cite{sha14}, the best bound one could get using a Metropolis weights is $\mathcal{O}(n^2\log n)$, as in Theorem \ref{teo2} if $B=1$.

 We note, however, that the requirements of Theorem \ref{thm_lineal} are more stringent than those of Theorem \ref{teo2}. The network topology is fixed (i.e. a static graph) and all nodes need to know an upper bound $U$ on the total number of agents. This upper bound must be within a constant factor of the number of agents. 
		
		
		\section{Generalized Distributed non-Bayesian Learning}\label{general}
		In this section, we discuss a general class of distributed non-Bayesian algorithms. 
		First, we will motivate the choice of the update rules described in  
		Eq.~(\ref{non_bayes_distributed}) and Eq.~\eqref{linear_bayes}. 
		For simplicity of exposition, we will assume that the agents
		always obtain observations (i.e. $\beta^i_k=1$ in Eqs.~(\ref{non_bayes_distributed}) and~\eqref{linear_bayes} 
		for all $i$ and $k$).
		Then, we will provide a comparison between our algorithms and previously proposed algorithms within
		the generalized distributed non-Bayesian framework.
		
		Standard centralized Bayes' rule can be described as the solution of a constrained optimization problem~\cite{zel88,wal06,hil12}. The cost function to be minimized is composed of two terms: one being the Maximum Likelihood Estimation (MLE) of a state given the observed data and the other being a regularization function minimized by the current prior~\cite{wal06}, i.e.,
		\begin{align*}
		\mu_{k+1}(\theta) & = \argmin_{\pi \in \mathbb{P}\left( \Theta\right) } \Big\{D_{KL}\left(\pi \| \mu_k \right) - \mathbb{E}_{\pi}[ \log\left(\ell\left(s_{k+1}|\theta\right)\right) ]\Big\} \\
		& = \frac{\mu_{k}(\theta)\ell(s_{k+1}|\theta)}{\sum_{p=1}^m\mu_{k}\left(\theta_p\right)\ell(s_{k+1}|\theta_p)}
		\end{align*}
		where $s_{k+1}$ is the most recent observation, $\ell(\cdot|\theta)$ is the likelihood function for hypothesis $\theta$, $\mathbb{E}_{\pi}$ is the expected value with respect to the probability distribution $\pi$, 
		and $\mathbb{P}\left( \Theta\right)$ is the set of all probability distributions on the set $\Theta$.
		
We can modify the optimization problem associated with a Bayesian update to take into account the network structure. We change the KL divergence term from a single prior belief to a convex combination of the beliefs of an agent and its neighbors in the network. The corresponding optimization problem for agent $i$ is: 
		\begin{align*}
		\mu_{k+1}^i(\theta) & 
		= \argmin_{\pi \in \mathbb{P}\left(\Theta\right)} \left\{\sum_{j=1}^n [A_k]_{ij} D_{KL}(\pi \| \mu_k^j )  \right. \\
		& \left. \qquad \qquad \qquad \vphantom{\sum_j a_{i,j} D_{KL}\left(\pi \| \mu_{j,k} \right)} 
		- \mathbb{E}_{\pi}[ \log\left(\ell^i(s_{k+1}^i\right|\theta))  ]\right\} \\
		& = \frac{\prod_{j=1}^n\mu_k^j(\theta)^{[A_k]_{ij}}\ell^i(s_{k+1}^i|\theta)}{\sum_{p=1}^m \prod_{j=1}^n \mu_k^j(\theta_p)^{[A_k]_{ij}}\ell^i(s_{k+1}^i|\theta_p)}.
		\end{align*}
		Observe that the solution of this optimization problem is precisely the proposed update rule in Eq.~\eqref{non_bayes_distributed}. 
		
		Opinion pooling or opinion aggregation has been studied before in~\cite{gen86,gil93,coo90,deg74}. 
		It is considered a traditional problem in economics, where several experts have beliefs 
		about a hypothesis and one needs to aggregate their beliefs into a single probability distribution. 
		Different opinion aggregation functions result from using different divergence metric for probability distributions 
		(see~\cite{gar04}). Similarly, different opinion pool operators define different non-Bayesian distributed learning rules. 
		A general form of opinion pooling was introduced in~\cite{gil93}, termed \textit{g-Quasi-Linear Opinion pools} (g-QLOP), 
		defined as follows:
		\begin{align*}
		\tau^{A_k}_g\left(\hdots,\mu_k^j (\theta),\hdots\right) & = \frac{g^{-1}\left(\sum_{j=1}^n {[A_k]_{ij}}g(\mu_k^j(\theta))\right)}{\sum_{p=1}^m g^{-1}\left(\sum_{j=1}^n {[A_k]_{ij}}g(\mu_k^j\left(\theta_p\right))\right)}
		\end{align*}
		with $\tau^A_g: \prod_{i=1}^n \mathbb{P}\left( \Theta\right)  \to \mathbb{P}\left(\Theta \right) $. The \textit{g-QLOP} corresponds to weighted arithmetic averages when $g(x) = x$ and 
		to weighted geometric averages when $g(x) = \log x $.
		
		The update rules studied in this paper can be seen as a two-step procedure. First, the beliefs of the neighbors are combined according to an opinion aggregation function. Second, the resulting aggregate distribution is updated using Bayes' rule. The proposed update rule, see Eq. (\ref{non_bayes_distributed}), uses the Logarithmic Opinion Pool, where
		\begin{align*}
		\tau^{A_k}_{\log x}\left(\hdots,\mu_k^j (\theta),\hdots\right) & = \frac{\prod_{j=1}^n \mu_k^j(\theta)^{[A_k]_{ij}}}{\sum_{p=1}^m\prod_{j=1}^n \mu_k^j\left(\theta_p\right)^{[A_k]_{ij}}}
		\end{align*}
		thus
		\begin{align*}
		\mu_{k+1}^i(\theta)  & = \frac{\tau^{A_k}_{\log x}\left(\hdots,\mu_k^j (\theta),\hdots\right) \ell^i(s_{k+1}^i|\theta)}{\sum_{p=1}^m \tau^{A_k}_{\log x}\left(\hdots,\mu_k^j (\theta_p),\hdots\right) \ell^i(s_{k+1}^i|\theta_p)}.
		\end{align*}
		
Logarithmic Pools are externally Bayesian \cite{gen86,mad64}, i.e. the order of aggregation of beliefs and the of new evidence does not influence the update rule. That is, from a learning point of view, if the function is Externally Bayesian, we can interchange the innovation and diffusion steps. The order in which we aggregate opinions and make the Bayesian update does not change the update rule. The next proposition shows that the update rule in Eq.~\eqref{non_bayes_distributed} is externally Bayesian.
		\begin{proposition} Assume that $\beta^i_k=1$ for all $i$ and $k$ in the update rule Eq.~\eqref{non_bayes_distributed}.
			Then, this rule is externally Bayesian, i.e. Eq.~\eqref{non_bayes_distributed} is equivalent to:
			\begin{align*}
			\mu_{k+1}^i(\theta)  & = \tau^{A_k}_{\log x}\left(\hdots,\frac{\mu_{k}^j(\theta)\ell^i(s^i_{k+1}|\theta)}{\sum_{p=1}^m\mu_{k}^j\left(\theta_p\right)\ell^i(s_{k+1}^i|\theta_p)} ,\hdots\right). 
			\end{align*}
		\end{proposition}
		\begin{IEEEproof}
			First generate a posterior taking as prior each of the opinions in the neighbor set:
			\begin{align*}
			\mu^i_{j,k+1}(\theta) 
			& = \frac{\mu_{k}^j(\theta)\ell^i(s^i_{k+1}|\theta)}
			{\sum_{p=1}^m \mu_{k}^j\left(\theta_p\right)\ell^i\left(s^i_{k+1}|\theta_p\right) }.
			\end{align*}
			Then combine the resulting $\mu^i_{j,k+1}$ into a new posterior, denoted by $\tilde \mu_{k+1}^i(\theta)$,
			as follows:
			\begin{align*}
			\tilde \mu_{k+1}^i(\theta) 
			& = \tau^{A_k}_{\log x}\left(\hdots,\mu_{j,k+1}^i (\theta),\hdots\right) \\
			& = \frac{\prod_{j=1}^n \mu^i_{j,k+1}(\theta)^{[ A_k]_{ij}}  }
			{\sum_{p=1}^m \prod_{j=1}^n \mu^i_{j,k+1}\left(\theta_p\right)^{[ A_k]_{ij}} }. 
			\end{align*}
			Substitute the expressions for $\mu^i_{j,k+1}(\theta)$ in the preceding relation to obtain
			\begin{align*}
			\tilde \mu_{k+1}^i(\theta) 
			& =  \frac{\prod_{j=1}^n \left(\frac{\mu_{k}^j(\theta)\ell^i(s^i_{k+1}|\theta)}
				{\sum_{q=1}^m \mu_{k}^j\left(\theta_q\right)\ell^i\left(s^i_{k+1}|\theta_q\right) }\right)^{[ A_k]_{ij}}  }
			{\sum_{p=1}^m \prod_{j=1}^n \left(\frac{\mu_{k}^j\left(\theta_p\right)\ell^i\left(s^i_{k+1}|\theta_p\right)}
				{\sum_{q=1}^m \mu_{k}^j\left(\theta_q\right)\ell^i\left(s^i_{k+1}|\theta_q\right) }\right)^{[ A_k]_{ij}} } \\
			& = 
			\frac{ \prod_{j =1}^n \mu_{k}^j(\theta)^{[ A_k]_{ij}}\ell^i(s^i_{k+1}|\theta)}
			{\sum_{p=1}^m \prod_{j =1}^n \mu_{k}^j\left(\theta_p\right)^{[ A_k]_{ij}}\ell^i\left(s^i_{k+1}|\theta_p\right)}
			\end{align*}
			where the last equality is obtained by noting that the term 
			$\prod_{j=1}^n ( \sum_{q=1}^m \mu_{k}^j\left(\theta_q\right)\ell^i(s^i_{k+1}|\theta_q) )^{[ A_k]_{ij} }$ cancels out from the numerator and the denominator.
			The last relation is the same as Eq.~(\ref{non_bayes_distributed}), so that 
			$\tilde \mu_{k+1}^i(\theta) =\mu_{k+1}^i(\theta)$.
		\end{IEEEproof}
		
		Consider now a Linear Opinion pool, where
		\begin{align*}
		\tau_x^{A_k}\left(\hdots,\mu_k^j (\theta),\hdots\right) & = \sum_{j=1}^n [A_k]_{ij} \mu_k^j(\theta).
		\end{align*}
		If the opinion aggregation is done first, as studied in \cite{qip15}, then
		the resulting update rule is
		\begin{align*}
		\mu_{k+1}^i(\theta) =\frac{\sum_{j=1}^n [A_k]_{ij}\mu_k^j(\theta)\ell^i(s_{k+1}^i|\theta)}{\sum_{p=1}^m \sum_{j=1}^n [A_k]_{ij}\mu_k^j(\theta_p) \ell^i(s_{k+1}^i|\theta_p)}.
		\end{align*}
		On the other hand, if the Bayesian update is done first,
		then the resulting update rule is
		\begin{align}\label{ATC_lin}
		\mu_{k+1}^i(\theta) = \sum_{j=1}^n [A_k]_{ij} \frac{\mu_k^j(\theta)\ell^j(s_{k+1}^j|\theta)}{\sum_{p=1}^m \mu_k^j(\theta_p) \ell^j(s_{k+1}^j|\theta_p)}.
		\end{align}
		
		The Linear Pool-based update rule is similar to the update rule proposed in~\cite{jad12}.  The authors in \cite{jad12} proposed the following rule
		\begin{align*}
		\mu_{k+1}^i(\theta) & = \tau_{x}^{A}\left(\hdots,\frac{\mu_{k}^i(\theta)\ell^i(s^i_{k+1}|\theta)}{\sum_{p=1}^m\mu_{k}^j\left(\theta_p\right)\ell^i(s_{k+1}^i|\theta_p)},\hdots\right)
		\end{align*}
		where opinion aggregation with linear functions is performed locally with priors from the neighbors. 
		The main difference is that in Eq.~\eqref{ATC_lin}, 
		a convex combination of the posteriors received from the neighbor set is used to generate the new individual posterior, while  
		in~\cite{jad12} the update rule is a convex combination of the individual posterior and the neighbors' priors.

In \cite{sha13}, the authors considered the case where the randomized gossip algorithm defines the communication structure. The update protocol is based on a distributed version of the Nesterov's dual averaging with stochastic gradients corresponding to the log-likelihood models given a set of observations. In this case, the agents exchange the likelihoods of the current observations instead of the beliefs. Thus, the consensus step is performed as a geometric aggregation of the likelihoods, and the resulting update rule can be described as
		\begin{align}\label{share_likeli}
		\mu_{k+1}^i(\theta) & = \frac{\mu_{k}^i(\theta)\tau_{\log x}^{W_k}\left(\hdots,\ell^j(s_{k+1}^j|\theta),\hdots\right)}{\sum_{p=1}^{m}\mu_{k}^i\left(\theta_p\right)\tau_{\log x}^{W_k}\left(\hdots,\ell^j(s_{k+1}^j|\theta_p),\hdots\right)}
		\end{align}
		where $W_k$ is the communication matrix coming from the gossip protocol.
		
		The idea of communicating aggregated versions likelihoods instead of beliefs was previously studied in the context of distributed estimation in sensor networks \cite{sal06}. Approaching the problem from the point of view of the Belief Propagation algorithm, resulted in an update rule in the form of Eq.~\eqref{share_likeli}. In \cite{sal06}, the authors showed convergence results for primitive, rings, tree, random graphs and other extensions to the original belief propagation algorithm. Similarly, in \cite{lal14b}, the authors propose an update rule where every agent performs local Bayesian updates before aggregating their beliefs using geometric averages, i.e.
		\begin{align*}
		\mu_{k+1}^i(\theta)  & = \tau^{A}_{\log x}\left(\hdots,\frac{\mu_{k}^j(\theta)\ell^j(s^j_{k+1}|\theta)}{\sum_{p=1}^m\mu_{k}^j\left(\theta_p\right)\ell^j(s_{k+1}^j|\theta_p)} ,\hdots\right).
		\end{align*}
		
		Convergence results for fixed communications matrices are provided, as well as asymptotic characterizations of the rates of convergence. Later in \cite{lal14b}, the authors extended the characterization of the rate of convergence to large deviation theory, providing a statement about the existence of a random time after which the beliefs will decrease exponentially.
		
		\section{Consistency of the Learning Rule}\label{consistency} 
		
		This section provides the proof for Theorem~\ref{proof_consistency}.
		We begin with a sequence of auxiliary lemmas. First, we recall few results from \cite{ned13} 
		about the convergence of a product of doubly stochastic matrices. 
		
		\begin{lemma}\label{lemma_angelia}
			\cite{ned13,ned09} 
			Under Assumption~\ref{assum:graph} on a matrix sequence $\{A_k\}$,
			we have
			\begin{align*}
			\left|  [A_{k:t}]_{ij} - \frac{1}{n}\right|  \leq \sqrt{2} \lambda^{k-t} \ \ \ \ \forall \ k \geq t \geq 0
			\end{align*}
			where $\lambda \in \left(0,1 \right) $ satisfies the relations described in Theorem~\ref{teo2}.
		\end{lemma}
		
		\begin{IEEEproof} 
			The proof may be found in~\cite{ned13}, with the exception of the bounds on $\lambda$ for the lazy Metropolis chains which may be found in \cite{ols14}.
		\end{IEEEproof}
		
		Next, we present a result regarding the weighted average of random variables with a finite variance.
		
		\begin{lemma}\label{waverage} 
			Assume that the graph sequence $\{\mathcal{G}_k\}$ satisfies Assumption \ref{assum:graph}. Also, let
			Assumption~\ref{assum:init} hold.
			Then, for $\theta_v \notin \Theta^*$ and $\theta_w \in \hat{\Theta}^*$, 
			\begin{align}\label{lemma:wave}
			\lim_{k \to \infty} \left(
			\frac{1}{k}\sum_{t=1}^{{ k}} A_{k:t} \mathcal{L}_{t}^{\theta_v,\theta_w} + 
			\frac{1}{n}\boldsymbol{1}_n\boldsymbol{1}_n' \boldsymbol{H}(\theta_v,\theta_w)  \right)
			= 0 \ a.s.
			\end{align}
			where $\mathcal{L}_t^{\theta_v,\theta_w}$ is a random vector with coordinates given by 
			\[ [ \mathcal{L}_t^{\theta_v,\theta_w} ]_i=  
			\beta^i_{t-1} \log \frac{\ell^i(S_{t}^i|\theta_v)}{\ell^i(S_{t}^i|\theta_w)} \qquad\forall i=1,\ldots,n\]
			while the vector $\boldsymbol{H}(\theta_v,\theta_w)$ has coordinates given by
			\begin{align*}
			H^i(\theta_v,\theta_w) = 
			q^i \left( D_{KL}(f^i \| \ell^i(\cdot|\theta_v)) - D_{KL}\left(f^i \|\ell^i\left(\cdot|\theta_w\right)\right)\right). 
			\end{align*}
		\end{lemma}
		\begin{IEEEproof}
			Adding and subtracting 
			$\frac{1}{k}{  \sum_{t=1}^k \frac{1}{n}\boldsymbol{1}_n\boldsymbol{1}_n' \mathcal{L}_{t}^{\theta_v,\theta_w}  }$ 
			to the expression under the limit in Eq.~(\ref{lemma:wave}) yields
			\begin{align}\label{eq:rel} 
			&\frac{1}{k}\sum_{t=1}^{{ k}} A_{k:t} 
			\mathcal{L}_{t}^{\theta_v,\theta_w} + \frac{1}{k}{  \sum_{t=1}^k \frac{1}{n}\boldsymbol{1}_n\boldsymbol{1}_n' \boldsymbol{H}(\theta_v,\theta_w) } 
			= \nonumber \\
			& \ \ \ \ \ \ \ \ \ \ \ \ \frac{1}{k}\sum_{t=1}^{{ k}} \left( A_{k:t} -\frac{1}{n}\boldsymbol{1}_n\boldsymbol{1}_n' \right) 
			\mathcal{L}_{t}^{\theta_v,\theta_w}   \nonumber \\
			& \ \ \ \ \ \ \ \ \ \ \ \ + \frac{1}{k} {  \sum_{t=1}^k \frac{1}{n}\boldsymbol{1}_n\boldsymbol{1}_n'    \left( \mathcal{L}_{t}^{\theta_v,\theta_w} +\boldsymbol{H}(\theta_v,\theta_w)\right)  }.
			\end{align}
			By Lemma \ref{lemma_angelia}, 
			$\lim_{k \to \infty} A_{k:t} =\frac{1}{n}\boldsymbol{1}_n \boldsymbol{1}_n'$ for all $t\ge0$. 
			Moreover, by Assumption~\ref{assum:init}(b), we have that $\log\alpha \leq [ \mathcal{L}_t^{\theta_v,\theta_w} ]_i \leq \log \frac{1}{\alpha}$. Thus, the first term on the right hand side of Eq.~\eqref{eq:rel} goes to zero a.s. as we take the limit over $k \to \infty$. 
			
			Regarding the second term on the right side of  Eq.~\eqref{eq:rel}, by the definition of the KL divergence, 
			and the assumption of each $\beta_t^i$ being independent, we have that   
			\begin{align*}
			& \mathbb{E}\left[ \beta^i_{t-1} \log \frac{\ell^i(S_{t}^i|\theta_v)}{\ell^i(S_{t}^i|\theta_w)}\right] 
			= q^i\sum_{s \in \mathcal{S}^i} f^i(s)   \log \frac{\ell^i\left(s|\theta_v\right)}{\ell^i\left(s|\theta_w\right)} \\
			&  = q^i \sum_{s \in \mathcal{S}^i} f^i(s)\log 
			\left( \frac{\ell^i\left(s|\theta_v\right)}{\ell^i\left(s|\theta_w\right)}\frac{f^i(s)}{f^i(s)} \right) \nonumber \\
			&  = q^i\left(  \sum_{s \in \mathcal{S}^i} f^i(s)\log 
			\left( \frac{f^i(s)}{\ell^i(s|\theta_w)} \right) -  \sum_{s \in \mathcal{S}^i} f^i(s)\log 
			\left( \frac{f^i(s)}{\ell^i(s|\theta_v)} \right)  \right)   \\
			& = q^i\left( D_{KL}\left(f^i\|\ell^i\left(\cdot|\theta_w\right)\right) - D_{KL}\left(f^i \|\ell^i\left(\cdot|\theta_v\right)\right)\right)
			\end{align*} 
			or equivalently
			\[ \mathbb{E} [ \mathcal{L}_t^{\theta_v,\theta_w} ]= -\boldsymbol{H}(\theta_v,\theta_w).\]
			
			Kolmogorov's strong law of large numbers states that if $\left\{X_t\right\}$ is a sequence of 
			independent random variables with variances such that $\sum_{k=1}^{\infty} \frac{{\rm Var}\left(X_k\right)}{k^2} < \infty$,
			then ${\frac{1}{n}\sum_{k=1}^n X_k - \frac{1}{n}\sum_{k=1}^n\mathbb{E}[ X_k]\to 0}$ a.s.
			Let ${X_t=\frac{1}{n}\boldsymbol{1}_n' \mathcal{L}_t^{\theta_v,\theta_w }}$, then by Assumption ~\ref{assum:init}(b), it can be seen
			that $\sup_{t\ge0} {\rm Var}\left(X_t\right) < \infty $. The final result follows by Lemma~\ref{lemma_angelia} and Kolmogorov's strong law of large numbers.
		\end{IEEEproof}
		
		Lemma~\ref{waverage} provides the necessary results to complete the proof of Theorem~\ref{proof_consistency}.

		\begin{IEEEproof}
			(\textit{Theorem~\ref{proof_consistency}}) 
			Initially, lets define the following quantities:
			for all $i=1,\ldots,n$ and $k\ge0$,
			\begin{align}\label{eq:defphi}
			\varphi_{k}^i (\theta_v,\theta_w)\triangleq\log \frac{\mu_{k}^i(\theta_v)}{\mu_{k}^i(\theta_w)}
			\end{align}
			defined for any $\theta_v \notin \hat{\Theta}^*$ and $\theta_w\in\hat{\Theta}^*$.
			We also use these quantities later in the proof of Theorem~\ref{teo2}.
			
			Let agent $i$ be arbitrary and consider the update rule of Eq.~(\ref{non_bayes_distributed}). 
			We will show that $\mu_k^i(\theta_v) \rightarrow 0$ as
			$k \rightarrow \infty$ for all $i=1, \ldots, n$. Note that if $\theta_v \in \Theta^* \setminus \hat\Theta^*$, then as a consequence of 
			Assumption \ref{assum:init}(a) we have that $\mu_k^i(\theta_v)=0$ for all $i$ and large enough $k$. Thus, we consider the case when $\theta_v \notin \Theta^*$ in the remainder of
			this proof. 
			
			Using the definition of $\varphi_{k}^i (\theta_v,\theta_w)$, it follows from Eq.~\eqref{non_bayes_distributed} that
			\begin{align*}
			\varphi_{k+1}^i(\theta_v,\theta_w)
			& = \log \frac{\mu_{k+1}^i(\theta_v)}{\mu_{k+1}^i(\theta_w)} \\
			& = 
			\log \frac{ \prod_{j=1}^n \mu_k^j(\theta_v)^{[A_k]_{ij}}\ell^i(S_{k+1}^i|\theta_v)^{\beta^i_k} }
			{\prod_{j=1}^n \mu_k^j(\theta_w)^{[A_k]_{ij}}\ell^i(S_{k+1}^i|\theta_w)^{\beta^i_k}} \\
			& = \sum_{j=1}^n [A_k]_{ij}\varphi_k^j\left(\theta_v, \theta_w\right) 
			+ \beta^i_k\log \frac{\ell^i(S_{k+1}^i|\theta_v)}{\ell^i(S_{k+1}^i|\theta_w)}. 
			\end{align*}
			
			Stacking up the values $\varphi_{k+1}^i(\theta_v,\theta_w)$ for $i=1, \ldots, n,$ 
			into a single vector $\boldsymbol{\varphi}_{k+1}(\theta_v,\theta_w)$, we can compactly write the preceding relations, as follows:
			\begin{align}\label{e6}
			\boldsymbol{\varphi}_{k+1}(\theta_v,\theta_w)
			& = A_k \boldsymbol{\varphi}_k(\theta_v,\theta_w) + \mathcal{L}_{k+1}^{\theta_v,\theta_w} 
			\end{align} 
			where $\mathcal{L}_{k+1}^{\theta_v,\theta_w}$ is defined in the statement of Lemma \ref{waverage}. Now, the relation in Eq.~\eqref{e6} implies that for all $k\ge0,$
			\begin{align}\label{e7}
			\boldsymbol{\varphi}_{k+1}(\theta_v,\theta_w)
			& = A_{k:0} \boldsymbol{\varphi}_0(\theta_v,\theta_w) 
			+ \sum_{t=1}^{k} A_{k:t} \mathcal{L}_{t}^{\theta_v,\theta_w} + \mathcal{L}_{k+1}^{\theta_v,\theta_w}.  
			\end{align} 
			
			The, if we add and subtract $\sum_{t=1}^{k}\frac{1}{n}\boldsymbol{1}_n\boldsymbol{1}_n' \boldsymbol{H}(\theta_v,\theta_w)$ in Eq.~\eqref{e7}, where $\boldsymbol{H}(\theta_v,\theta_w)$ is as in Lemma~\ref{waverage}, it follows that 
				\begin{align*}
				& \boldsymbol{\varphi}_{k+1}(\theta_v,\theta_w)
				=A_{k:0} \boldsymbol{\varphi}_0(\theta_v,\theta_w) 
				-  \frac{k}{n} \sum\limits_{i=1}^{n}H^i(\theta_v,\theta_w) \boldsymbol{1}_n  \\
				& \qquad
				+ \sum_{t=1}^{{ k}} \left( A_{k:t} \mathcal{L}_t^{\theta_v,\theta_w} 
				+ \frac{1}{n}\boldsymbol{1}_n\boldsymbol{1}_n' \boldsymbol{H}(\theta_v,\theta_w) \right) + \mathcal{L}_{k+1}^{\theta_v,\theta_w}  .
				\end{align*}

			By the definition of group confidence (cf.~Definition~\ref{confi}), we have
			\begin{align}\label{eq:hc}\sum_{i=1}^{n}H^i(\theta_v,\theta_w) 
			=\mathsf{C}_{\boldsymbol{q}}(\theta_w) -\mathsf{C}_{\boldsymbol{q}}(\theta_v)
			=\mathsf{C}^*_{\boldsymbol{q}}-\mathsf{C}_{\boldsymbol{q}}(\theta_v)
			\end{align}
			where the last equality follows from $\theta_w\in  \hat \Theta^*$ and 
			the definition of the optimal value $\mathsf{C}_{\boldsymbol{q}}^*$ (Definition~\ref{def:opt}).
			Therefore,
				\begin{align*}
				&\boldsymbol{\varphi}_{k+1}(\theta_v,\theta_w) 
				 =  A_{k:0} \boldsymbol{\varphi}_0(\theta_v,\theta_w)  
				-  \frac{k}{n} \left(\mathsf{C}^*_{\boldsymbol{q}}
				- \mathsf{C}_{\boldsymbol{q}}(\theta_v)\right) \boldsymbol{1}_n \\
				& \qquad
				+ \sum_{t=1}^{{ k}} \left( A_{k:t} \mathcal{L}_t^{\theta_v,\theta_w} +\frac{1}{n}\boldsymbol{1}_n\boldsymbol{1}_n' \boldsymbol{H}(\theta_v,\theta_w) \right) + \mathcal{L}_{k+1}^{\theta_v,\theta_w}    .
				\end{align*}

			By dividing both sides of the preceding equation with $k$ and taking the limit as $k$ goes to infinity, almost surely we have
			\begin{align}\label{eq8}
			&\lim_{k \to \infty} \frac{1}{k} \boldsymbol{\varphi}_{k+1}(\theta_v,\theta_w)
			=  \lim_{k \to \infty} \frac{1}{k} A_{k:0} \boldsymbol{\varphi}_0(\theta_v,\theta_w) \nonumber\\
			&
			\ \ \ \ + \lim_{k \to \infty} \frac{1}{k} 
			\sum_{t=1}^{{ k}} \left(A_{k:t} \mathcal{L}_t^{\theta_v,\theta_w} +\frac{1}{n}\boldsymbol{1}_n\boldsymbol{1}_n' \boldsymbol{H}(\theta_v,\theta_w)\right)\cr
			& \ \ \ \ + \lim_{k \to \infty}\frac{1}{k} \mathcal{L}_{k+1}^{\theta_v,\theta_w}
			- \frac{1}{n}\left(\mathsf{C}_{\boldsymbol{q}}^*-\mathsf{C}_{\boldsymbol{q}}(\theta_v)\right) \boldsymbol{1}_n.
			\end{align}
			
			The limit on the left hand side of Eq.~\eqref{eq8} is justified since all the limits on the right-hand side exist. Specifically,
			the first term of the right-hand side of Eq.~\eqref{eq8} converges to zero deterministically. 
			The second term converges to zero almost surely by Lemma \ref{waverage}, while 
			the third term goes to zero since $\mathcal{L}_t^{\theta_v,\theta_w}$ is bounded almost surely (cf.\ Assumption~\ref{assum:init}(b)). 
			
			Consequently, 
			\begin{align*}
			\lim_{k \to \infty} \frac{1}{k}\boldsymbol{\varphi}_{k+1}(\theta_v,\theta_w) 
			= 
			- \frac{1}{n}\left(\mathsf{C}^*_{\boldsymbol{q}} 
			- \mathsf{C}_{\boldsymbol{q}}(\theta_v)\right) \boldsymbol{1}_n \ \ \text{a.s.}
			\end{align*} 
			
			Since $\mathsf{C}^*_{\boldsymbol{q}}$ is the maximum value and $\theta_v\not\in \Theta^*$, it follows that 
			$\mathsf{C}^*_{\boldsymbol{q}} -\mathsf{C}_{\boldsymbol{q}}(\theta_v)>0$, implying that 
			$\boldsymbol{\varphi}_k(\theta_v,\theta_w)\to-\infty$ almost surely.
			Also, by $ \mu^i_k(\theta_v)  \leq \exp\left(\varphi^i_k(\theta_v,\theta_w)\right)
			$ for all $i$, we have  $\mu^i_k(\theta_v)\to 0$ \it{a.s}. 
		\end{IEEEproof}

One specific instance of our setup is when there exists a unique hypothesis that matches the distribution of the observations of all agents. This case relates to the previously proposed approaches for distributed learning. Specifically, in~\cite{lal14b,sha13,jad12}, the authors assume that there is a ``true state" of the world, i.e., there is a unique hypothesis such that the distance between such hypothesis and the true distribution of the data is zero for all agents. This case could be expressed as a consequence Theorem~\ref{proof_consistency}, as follows:
		\begin{corollary}
			Under assumptions of Theorem~\ref{proof_consistency},
			if there is a unique hypothesis $\theta^*$ with $C_q^*=0$, then
			\begin{align*}
			\lim\limits_{k \to \infty } \mu_k^i(\theta^*) = 1 \ \ \  a.s.  \ \ \ \forall i \in V.
			\end{align*}
		\end{corollary}
		\begin{IEEEproof}
			By Theorem~\ref{proof_consistency} for every $\theta \neq \theta^*$ we have that 
			$		\lim\limits_{k \to \infty } \mu_k^i(\theta) = 0$ {\it a.s.}
		\end{IEEEproof}
		
\begin{figure}[H]
			\centering
			\includegraphics[width=0.4\textwidth]{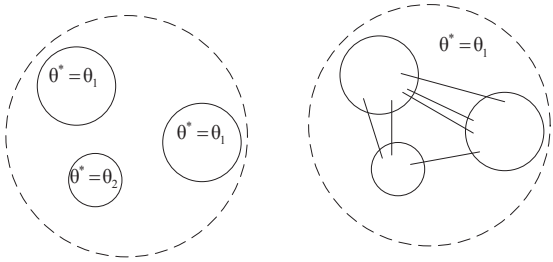}
			\caption{Conflicting social groups interacting. Initially on the left, there are three isolated social clicks, each with a different optimal hypothesis. Once such groups interact (on the right), others might influence the local decision and a click changes its beliefs to the optimal with respect to the complete set of agents. In this case, one of the groups was convinced that $\theta_1$ was a better solution than $\theta_2$.}
			\label{clicks}
		\end{figure}

In general, one can consider several closed social \textit{cliques} where the same hypothesis can represent different distributions for different groups. For example, in a social network, what one community might consider as a good hypothesis, need not be good for other communities. Each disconnected social \textit{clique} could have a different optimal hypothesis, even if all observations come from the same distribution, see Figure \ref{clicks}. If such social clicks interact, Theorem \ref{proof_consistency} provides the conditions for which all agents will agree on the a hypothesis that is the closest to the best one considering the models of all agents in the network and not only those in a specific \textit{clique}. 
   
		The previous statement is formally stated in the next corollary.
		\begin{corollary}\label{cor_clicks}
			Let the agent set $V$ be partitioned into $\hat p$ disjoint sets $V_j,j=1,\ldots,\hat p$.	Under assumptions of Theorem~\ref{proof_consistency} where each agent updates 
			its beliefs according to Eq.~\eqref{non_bayes_distributed}, 
			if there exists a hypothesis $\theta^*$ such that
			\begin{align*}
			\sum_{j=1}^{\hat p} C_q^{V_j}(\theta^*) & > \max_{\theta \neq \theta^*} \sum_{j=1}^{\hat p}  C_q^{V_j}(\theta),
			\end{align*}
			then $\lim_{k \to \infty} \mu_k^i(\theta^*)  = 1 $ \text{a.s.  for all $i$}. 
		\end{corollary}
		\begin{IEEEproof}
			If the hypothesis $\theta^*$ exists, then
			the group confidence on $\theta^*$ is larger than the group confidence for any other hypothesis. 
			Thus, $\hat{\Theta}^*=\{\theta^* \}$ and the result follows by Theorem \ref{proof_consistency}.
		\end{IEEEproof}
		
		\section{Rate of Convergence for Time-Varying Graphs}\label{rates}
		In this section, we prove Theorem \ref{teo2}, which provides an explicit rate of convergence for the learning process. 
		
		The next lemma is an extension of Lemma 2 in \cite{sha14} to the case of time-varying graphs. It provides a technical result that will help us later on the computation of the non-asymptotic convergence rate.
		\begin{lemma}\label{shahin}
			Let Assumption \ref{assum:graph} hold for a matrix sequence $\{A_k\}$. Then for all $i$, 
			\begin{align*}
			\sum\limits_{t=1}^{k} \sum\limits_{j=1}^{n} \left| [A_{k:t}]_{ij} - \frac{1}{n} \right| & \leq \frac{4 \log n}{1-\lambda}
			\end{align*}
			where $\lambda = 1-\eta / 4n^2$, and if every $A_k$ is a lazy Metropolis matrix then $\lambda = 1- 1 / \mathcal{O}(n^2)$.
		\end{lemma}
		\begin{IEEEproof}
			In \cite{sha14}, the authors assume the weight matrix is static and diagonalizable, then they use the following inequality from \cite{ros95}:
			\begin{align*}
				\|\boldsymbol{e}_j'A^k - \pi'\|_1 & \leq n \lambda_{\text{max}}(A)^k
			\end{align*}
			where $\boldsymbol{e}_j$ is a vector with its $j$-th entry equal to one and zero otherwise, $\pi$ is the stationary distribution of the Markov Chain with transition matrix $A$ and $\lambda_{\text{max}}(A)$ is the second largest eigenvalue of the matrix $A$.

		For time-varying graphs one can use the inequality in Lemma \ref{lemma_angelia} instead. The reminder of the proof remains the same as in \cite{sha14}.
		\end{IEEEproof}
		
		Before proving Theorem \ref{teo2}, we will provide an auxiliary result regarding bounds on the expectation of the random variables $\varphi_k^i(\theta_v,\theta_w)$ as defined in Eq.~\eqref{eq:defphi}. 
		
		\begin{lemma}\label{bound_vhi}
			Consider $\varphi^i_{k}(\theta_v,\theta_w)$ as defined in Eq.~\eqref{eq:defphi}, 
			with $\theta_w\in\hat\Theta^*$.
			Then, for any $\theta_v\notin \Theta^*$  we have
			\begin{align*}
			\mathbb{E}[ \varphi^i_{k}\left(\theta_v,\theta_w \right)]& \leq
			\gamma_1^i - k \gamma_2\quad\hbox{for all $i$ and $k\ge 0$}
			\end{align*}
			with $\gamma_1^i$ and $\gamma_2$ as defined in Theorem \ref{teo2}.
		\end{lemma}
		\begin{IEEEproof}
			Taking the expected value in Eq.~\eqref{e7} we can see that for all $k\ge0$,
			\begin{align*}
			\mathbb{E}[ \varphi_{k+1}^i(\theta_v,\theta_w)]
			& = \sum_{j=1}^n [ A_{k:0}]_{ij} \varphi_0^j(\theta_v,\theta_w)  \\
			&  -\sum_{t=1}^{k} \sum\limits_{j=1}^{n} [ A_{k:t}]_{ij} H^j(\theta_v,\theta_w) 
			- H^i(\theta_v,\theta_w).
			\end{align*}
			
			By adding and subtracting $\sum\limits_{t=1}^{k+1} \sum\limits_{j=1}^{n} \frac{1}{n} H^j(\theta_v,\theta_w)$, we obtain
			{\small
			\begin{align}\label{exp_1}
			& \mathbb{E}[ \varphi_{k+1}^i(\theta_v,\theta_w)]
			 = \sum_{j=1}^n [ A_{k:0}]_{ij} \varphi_0^j(\theta_v,\theta_w)  \nonumber\\
			& \ - \sum_{t=1}^{k} \sum\limits_{j=1}^{n}\left(  [ A_{k:t}]_{ij} - \frac{1}{n} \right) H^j(\theta_v,\theta_w)  - \frac{k+1}{n} \sum\limits_{j=1}^{n} H^j(\theta_v,\theta_w) \nonumber\\
			& \qquad 
			- \left( H^i(\theta_v,\theta_w) - \frac{1}{n}\sum\limits_{j=1}^{n} H^j(\theta_v,\theta_w)\right).
			\end{align}
			}
			
			For the first term in Eq. \eqref{exp_1}, since $A_{k:0}$ is stochastic matrix, we have that
			\begin{align*}
			\sum_{j=1}^n [ A_{k:0}]_{ij} \varphi_0^j(\theta_v,\theta_w) & \leq  \max_i\log \frac{\mu_0^i(\theta_v)}{\mu_0^i(\theta_w)}  .
			\end{align*}
			
			The second term in Eq. \eqref{exp_1} can be bounded using Lemma~\ref{shahin}, thus
			\begin{align*}
			\sum_{t=1}^{k} \sum\limits_{j=1}^{n}\left(  [ A_{k:t}]_{ij} - \frac{1}{n} \right) H^j(\theta_v,\theta_w) & \leq \frac{4 \log n}{1-\lambda} \log\frac{1}{\alpha}
			\end{align*}
			since $\log \alpha \leq H^j(\theta_v,\theta_w) \leq \log\frac{1}{\alpha} $.
			
			The last term in Eq. \eqref{exp_1} is bounded as
			\begin{align*}
			\frac{1}{n}\sum\limits_{j=1}^{n} \left( H^i(\theta_v,\theta_w) - H^j(\theta_v,\theta_w)\right)&\leq  2\log \frac{1}{\alpha} \\
			&\leq\frac{8\log n}{1-\lambda}\log \frac{1}{\alpha}
			\end{align*}
			where the last inequality follows from $2 \leq 8\log n$ for $n\geq 2$ and $1-\lambda <1$.
			
			Finally we have that
				\begin{align*}
				\mathbb{E}[ \varphi_{k+1}^i(\theta_v,\theta_w)]
				& \leq \max_i\log \frac{\mu_0^i(\theta_v)}{\mu_0^i(\theta_w)}  + \frac{12 \log n}{1-\lambda} \log\frac{1}{\alpha} \\
				& \qquad \qquad \qquad - \frac{k+1}{n} \sum\limits_{j=1}^{n} H^j(\theta_v,\theta_w)
				\end{align*}
				from which the desired result follows by using the definitions of 
				$\gamma_1^i$, $\gamma_2$, $H^j(\theta_v,\theta_w)$ and taking the appropriate maximum values over $\theta_v$ and $\theta_w$ on the right hand side of the preceding inequality.
		\end{IEEEproof}
		
In the proof of Theorem \ref{teo2} we will use McDiarmid's inequality~\cite{mcd89}, which provides bounds for the concentration of functions of random variables. This inequality allows us to show bounds on the probability that the beliefs exceed a given value $\epsilon$. For completeness, next we state the McDiarmid's inequality.
		
		\begin{theorem}\label{mcd}
			(McDiarmid's inequality \cite{mcd89})
			Let $X_1,\hdots,X_k$ 
			be a sequence of independent random variables with $X_t \in \mathcal{X}$ for $1 \leq t \leq k$. Further, let 
			$g: \mathcal{X}^k \to \mathbb{R}$ be a function of bounded differences, i.e.,
			for all $1 \leq t \leq k$,
			\begin{align*}
			\sup\limits_{x_t \in \mathcal{X}} g\left(\hdots,x_t,\hdots \right) 
			-\inf\limits_{x_t \in \mathcal{X}} g\left(\hdots,x_t,\hdots \right)  &\leq c_t
			\end{align*}
			then for any $\epsilon > 0$ and all $k\ge1$,
			\begin{align*}
			\mathbb{P}\left(g(\left\lbrace X_t\right\rbrace_{t=1}^k  ) 
			- \mathbb{E}[g(\left\lbrace X_t\right\rbrace_{t=1}^k ) ] \geq \epsilon \right) 	\leq \exp\left(-\frac{2\epsilon^2}{\sum_{t=1}^kc_t^2} \right). 
			\end{align*}
		\end{theorem}
		Now, we are ready to prove Theorem~\ref{teo2}.
		
		\begin{IEEEproof}
			(\textit{Theorem} \ref{teo2}) First, we will express the belief $\mu_{k+1}^i(\theta_v)$ in terms of the variable $\varphi_{k+1}^i(\theta_v,\theta_w)$. This will allow us to use 
			McDiarmid's inequality to obtain the concentration bounds.
			By the dynamics of the beliefs in Eq.~\eqref{non_bayes_distributed} and Assumption~\ref{assum:init}(a), since
			$\mu_{k}^i(\theta_w)\in(0,1]$ for $\theta_w\in\hat\Theta^*$,  we have
			\begin{align*}
			\mu_{k}^i(\theta_v) 
			&\leq \frac{\mu_{k}^i(\theta_v)}{\mu_{k}^i(\theta_w)} 
			=\exp\left( \varphi^i_{k}(\theta_v,\theta_w)\right).
			\end{align*}
			
			Therefore,
			\begin{align*}
			&\mathbb{P}\left(\mu_{k}^i(\theta_v) 
			\geq \exp\left( -\frac{k}{2}\gamma_2 + \gamma_1^i \right) \right) \\
			& \ \ \ \ \ \ \leq  \mathbb{P}\left(\exp\left( \varphi_{k}^i(\theta_v,\theta_w) \right)   
			\geq \exp\left( -\frac{k}{2}\gamma_2 + \gamma_1^i  \right) \right) \\
			& \ \ \ \ \ \ =  \mathbb{P}\left(\varphi_{k}^i(\theta_v,\theta_w) \geq   -\frac{k}{2}\gamma_2 + \gamma_1^i \right)  \\
			& \ \ \ \ \ \ \leq \mathbb{P}\left(\varphi_{k}^i(\theta_v,\theta_w) - \mathbb{E}[\varphi_{k}^i(\theta_v,\theta_w) ] \geq   \frac{k}{2}\gamma_2  \right)
			\end{align*}
			where the last inequality follows from Lemma \ref{bound_vhi}.
			
			We now view $\varphi_{k+1}^i(\theta_v,\theta_w)$ as a function of the random vectors 
			$\boldsymbol{S}_1,\ldots,\boldsymbol{S}_k$, (see Eq.~\eqref{e7}), where  $\boldsymbol{S}_t=(S_t^1,\ldots,S_t^n)$ for $t\geq 1$, and the random variable $S_{k+1}^i$. 
			Next, we will establish that this function has bounded differences 
			in order to apply McDiarmid's inequality. 
			
			For all $t$ with $1\le t\le k$ and $j$ with $1\leq j \leq n$, we have
				\begin{align*}
				&\max_{s_t^j \in \mathcal{S}^j} \varphi_{k+1}^i(\theta_v,\theta_w) - \min_{s_t^j \in \mathcal{S}^j} \varphi_{k+1}^i(\theta_v,\theta_w) \\
				&  \ \ =  \max_{s_t\in \mathcal{S}^j} [  A_{k:t}] _{ij} \log \frac{\ell^j(s_t|\theta_v ) }{\ell^j(s_t|\theta_w )}    -\min_{s_t \in \mathcal{S}^j}  [  A_{k:t}] _{ij} \log \frac{\ell^j(s_t|\theta_v ) }{\ell^j(s_t|\theta_w )}  \\
				& \ \ \leq  [  A_{k:t}] _{ij} \log \frac{1 }{\alpha} +  [  A_{k:t}] _{ij}\log \frac{1}{\alpha}  \\
				& \ \ = 2  [  A_{k:t}] _{ij}\log \frac{1}{\alpha}. 
				\end{align*}
			
				Similarly, from Eq.~\eqref{e7} we can see that
				\[\max_{s_{k+1}^j \in \mathcal{S}^j} \varphi_{k+1}^i(\theta_v,\theta_w) 
				- \min_{s_{k+1}^j \in \mathcal{S}^j} \varphi_{k+1}^i(\theta_v,\theta_w) \le 2 \log \frac{1}{\alpha}.\]
			
			It follows that  $\varphi_{k+1}^i(\theta_v,\theta_w)$ has bounded variations, with
				\begin{align*}
				&\sum_{t=1}^{k}\sum_{j=1}^n (2  [  A_{k:t}] _{ij}\log \frac{1}{\alpha})^2 + \left( 2 \log \frac{1}{\alpha}\right)^2 \\
			& 	\qquad \qquad= 4\left(\log \frac{1}{\alpha}\right) ^2\left( \sum_{t=1}^{k}\sum_{j=1}^n ([  A_{k:t}] _{ij})^2+1\right) \\
				& 	\qquad \qquad\leq 4\left(\log \frac{1}{\alpha}\right) ^2(k+1)
				\end{align*}
				where the last inequality follows from the fact that $A_{k:t}$ is row stochastic.

		Thus,
			\begin{align*}
			& \mathbb{P}\left(\varphi_{k}^i(\theta_v,\theta_w) - \mathbb{E}[\varphi_{k}^i(\theta_v,\theta_w) ] \geq   \frac{k}{2}\gamma_2  \right) \\
			& \qquad \qquad =  \exp\left( -\frac{2\left( \frac{1}{2}k\gamma_2 \right)^2}{ 4 k\left( \log \frac{1}{\alpha}\right)^2}  \right). 
			\end{align*}
			
			Therefore, for a given confidence level $\rho$, in order to have 
			$\mathbb{P}\left(\mu_{k}^i(\theta_v) \geq \exp\left( -\frac{1}{2}k\gamma_2 + \gamma_1^i \right) \right) 
			\leq \rho$ we require that
			\begin{align*}
			k & \geq \frac{1}{\gamma_2^2}8 \left( \log\alpha \right) ^2  \log\frac{1}{\rho}.
			\end{align*}
		\end{IEEEproof}
		
		\section{Accelerated Learning for Fixed Undirected Graphs}\label{linear_rates}
		
		In this section, we analyze the distributed learning algorithm of Eq.~\eqref{linear_bayes} and prove its non-asymptotic convergence rate. 
		First, we will state an enabling theorem presented in~\cite{ols14}, 
		which presents a distributed consensus protocol that achieves a consensus with a linear growth in the number of agents.
		
		\begin{theorem}\label{alex_linear}
			\cite{ols14} Suppose each node $i$ in a fixed undirected connected graph updates its variable $x^i_k$ 
			at each time instant $k\ge 2$ as follows:
			\begin{subequations}\label{cons_upd}
				\begin{align}
				y^i_{k+1} & = x^i_k +\frac{1}{2}\sum_{j\in N_i} \frac{x^j_k - x^i_k}{\max\left\{d^i+1,d^j+1\right\}}\\
				x^i_{k+1} & = y^i_{k+1} + \left(1- \frac{2}{9U+1}\right)\left(y^i_{k+1} - y^i_{k}\right)
				\end{align}
			\end{subequations}
			where $N_i$ is the set of neighbors of agent $i$ and $d^i $ is its corresponding degree. Then, if $U \geq n$ we have that
			\begin{align}\label{thm_linear}
			\| \boldsymbol{y}_k - \bar{x} \boldsymbol{1} \|_2^2 & \leq 2 \left(1 - \frac{1}{9U}\right)^{k-1} \| \boldsymbol{y}_1 - \bar{x} \boldsymbol{1} \|_2^2	\quad \forall k\ge 1
			\end{align}
			where $[\boldsymbol{y}_k]_i = y^i_k $ and $\bar{x} = \frac{1}{n} \sum\limits_{i=1}^{n} x^i_1$, and 
			the process is initialized with $y^i_1= x^i_1$.
		\end{theorem}
		
		Next, we define some quantities that we use in the analysis of Eq. \eqref{linear_bayes}. Define the matrix $B$ and a scalar $\sigma$, as follows:
		\begin{align}\label{b_matrix}
		B & = \begin{bmatrix} \left(1+\sigma\right)\bar{A} & - \sigma \bar{A} \\ I_n & \textbf{0} \end{bmatrix}
		\end{align}
		\begin{align}\label{a_param}
		\sigma & = 1- \frac{2}{9U+1}
		\end{align}
		where $I_n$ is the identity matrix and $\textbf{0}$ is the matrix with all entries equal to zero of the appropriate size and $\bar{A}$ is as defined in Assumption \ref{assum_linear}. 
		
		We have the following auxiliary result for the matrix $B$.
		\begin{lemma}\label{analog_lemma3}
			Consider the matrix $B$ and the parameter $\sigma$ as defined in Eqs.~\eqref{b_matrix} and~\eqref{a_param} respectively. Then
			\begin{align*}
			\left| [[I_n \ \boldsymbol{0}]B^{k} [I_n \ I_n]']_{ij} - \frac{1}{n}\right|  \leq \sqrt{2} \lambda^{k} \ \ \ \ \forall \ k \geq 2
			\end{align*}
			where $ \lambda = 1 - \frac{1}{18U}$.
		\end{lemma}
		\begin{IEEEproof}
			The linear time consensus algorithm described in Eq.~\eqref{cons_upd} can be expressed as
			\begin{align*}
			\boldsymbol{y}_{k+1} & = \bar{A} \boldsymbol{x}_k \\
			\boldsymbol{x}_{k+1} & = \boldsymbol{y}_{k+1} + \sigma\left(\boldsymbol{y}_{k+1} - \boldsymbol{y}_{k}\right)
			\end{align*}                                         
			which implies that
			$\boldsymbol{y}_{k+1}  = \bar{A}\left(\boldsymbol{y}_{k} + \sigma\left(\boldsymbol{y}_{k} - \boldsymbol{y}_{k-1}\right)\right)$ with $y^i_1= x^i_1$.
			Therefore
			\begin{align*}
			\left[  \begin{array}{c} 	
			\boldsymbol{y}_{k+1}  \\ 	\boldsymbol{y}_{k}  \end{array} \right] & = \begin{bmatrix} \left(1+\sigma\right)\bar{A} & - \sigma \bar{A} \\ I_n & \textbf{0} \end{bmatrix}  \left[  \begin{array}{c} 	\boldsymbol{y}_{k}  \\ 	\boldsymbol{y}_{k-1}  \end{array} \right] = B^k  \left[  \begin{array}{c} 	\boldsymbol{y}_{1} \\ 	\boldsymbol{y}_{0}  \end{array} \right]
			\end{align*}
			where we assumed that $\boldsymbol{y}_0 = \boldsymbol{y}_1$. Thus,
			\begin{align*}
			\boldsymbol{y}_{k+1} & = [I_n \ \boldsymbol{0}]B^{k} [I_n  \ I_n]'\boldsymbol{y}_{1}.
			\end{align*}
			
			By substituting the previous relation into Eq.~\eqref{thm_linear} and using $\boldsymbol{x}_1=\boldsymbol{y}_1$, we obtain
			\begin{align}\label{square_bound}
			\|{ [I_n \ \boldsymbol{0}]B^{k} [I_n  \ I_n ]'\boldsymbol{y}_1 - \left( \frac{1}{n}\sum\limits_{i=1}^{n} y^i_1\right)  \boldsymbol{1}_n} \|_2^2 
			\nonumber \\
			\qquad \leq 2 \left(1 - \frac{1}{9U}\right)^{k} \| \boldsymbol{y}_1 - \frac{1}{n}\sum\limits_{i=1}^{n} y^i_1 \boldsymbol{1}_n \|_2^2
			\end{align}
			which implies that
			\begin{align*}
			&\max_{i} \left| [[I_n \ \boldsymbol{0}]B^{k} [I_n  \ I_n]'\boldsymbol{y}_1]_i - \frac{1}{n} \sum_{i=1}^n y_1^i\right| 
			\\ &\qquad \qquad \leq \sqrt{2} \left( \sqrt{1 - \frac{1}{9U}}\right)^{k} 
			\| \boldsymbol{y}_1 - \frac{1}{n}\sum\limits_{i=1}^{n} y^i_1 \boldsymbol{1}_n \|_2.
			\end{align*}
			
			The preceding relation holds for any $\boldsymbol{y}_1$. In particular, if we take  $\boldsymbol{y}_1=\boldsymbol{e}_j$, where $\boldsymbol{e}_j$ is a vector whose $j$-th entry is equal to one and zero otherwise, we conclude that for every $i$ and $j$,
			\begin{align*}
			\left| [[I_n \ \boldsymbol{0}]B^{k} [I_n  \ I_n]' ]_{ij} - \frac{1}{n}\right| 
			& \leq \sqrt{2} \left( 1 - \frac{1}{18U}\right)^{k}. 
			\end{align*}
			This follows from the inequality $\sqrt{1 - \beta} \leq 1 - \beta/2$ for all $\beta \in \left(0,1\right)$ and the fact that $\|\boldsymbol{e}_j - \frac{1}{n} \boldsymbol{1}_n\| \leq 1$.
		\end{IEEEproof}
		
		Now, we are ready to proof Theorem \ref{thm_lineal}.
		
		\begin{IEEEproof}
			(Theorem \ref{thm_lineal})	The proof is along the lines of the proof for Theorem~\ref{teo2}. From the definition of 
			$	\varphi_{k+1}^i(\theta_v,\theta_w) $ we have
			\begin{align*}
			&	\varphi_{k+1}^i(\theta_v,\theta_w)  = \log \frac{\mu_{k+1}^i(\theta_v)}{\mu_{k+1}^i(\theta_w)} \\
			& = \log \frac{ \frac{ \prod_{j=1}^n \mu_k^j (\theta_v)^{(1+\sigma) \overline{A}_{ij}}  \ell^i(S_{k+1}^i|\theta_v)^{\beta_k^i}  }{ \prod_{j=1}^n ( \mu_{k-1}^j(\theta_v) \ell^j(S_k^j | \theta_v)^{\beta_{k-1}^j} )^{\sigma \overline{A}_{ij}} }}{ \frac{ \prod_{j=1}^n \mu_k^j (\theta_w)^{(1+\sigma) \overline{A}_{ij}} \ell^i(S_{k+1}^i|\theta_w)^{\beta_k^i} }{ \prod_{j=1}^n( \mu_{k-1}^j(\theta_w) \ell^j(S_k^j | \theta_w)^{\beta_{k-1}^j})^{\sigma \overline{A}_{ij}}}}  \\ 
			& = \sum\limits_{j=1}^{n} 
			\left(1+\sigma\right) \bar{A}_{ij} \log \frac{\mu_{k}^j(\theta_v)}{\mu_{k}^j(\theta_w)} 
			- \sum\limits_{j=1}^{n} \sigma\bar{A}_{ij} \log \frac{\mu_{k-1}^j(\theta_v)}{\mu_{k-1}^j(\theta_w)} \\
			& \ \ \ \ \ 
			+ \beta^i_k\log \frac{\ell^i (S_{k+1}^i|\theta_v) }{\ell^i (S_{k+1}^i|\theta_w) } - \sum\limits_{j=1}^{n} \sigma\bar{A}_{ij} \beta^j_{k-1} \log \frac{\ell^j (S_{k}^j|\theta_v	) }{\ell^j (S_{k}^j|\theta_w) }  \\
			& =  \sum\limits_{j=1}^{n} \left(1+\sigma\right)\bar{A}_{ij} \varphi_{k}^j(\theta_v,\theta_w) 
			- \sum\limits_{j=1}^{n} \sigma\bar{A}_{ij} \varphi_{k-1}^j(\theta_v,\theta_w) 
			\\
			&\ \ \ \ \  + [\mathcal{L}_{k+1}^{\theta_v,\theta_w} ]_i - \sum\limits_{j=1}^{n} \sigma\bar{A}_{ij} [\mathcal{L}_{k}^{\theta_v,\theta_w} ]_j.
			\end{align*}
			Stacking the previous relation for all $i$ we obtain the following vector representation for the dynamics
			\begin{align}
			\boldsymbol{\varphi}_{k+1}(\theta_v,\theta_w) & = \left(1+\sigma\right)\bar{A} \boldsymbol{\varphi}_{k}(\theta_v,\theta_w) - \sigma\bar{A} \boldsymbol{\varphi}_{k-1}(\theta_v,\theta_w) \nonumber\\
			& \qquad  + \mathcal{L}_{k+1}^{\theta_v,\theta_w} - \sigma\bar{A}\mathcal{L}_{k}^{\theta_v,\theta_w}. 
			\end{align}
			Now, define the following auxiliary vector 
			\begin{align*}
			\boldsymbol{z}_{k+1}(\theta_v,\theta_w) & = \boldsymbol{\varphi}_{k}(\theta_v,\theta_w) +\mathcal{L}_{k+1}^{\theta_v,\theta_w}
			\end{align*} 
			where $\boldsymbol{z}_{0}(\theta_v,\theta_w) = 0$, since $\boldsymbol{\varphi}_{-1}(\theta_v,\theta_w)= 0$ by the assumption of uniform initial beliefs, and $\mathcal{L}_{0}^{\theta_v,\theta_w} = 0$ due to $\beta_{-1}=0$, in which case we can set $S^i_{0}$ to any value in $\mathcal{S}^i$.
			
			By writing the evolution for the augmented state $[ 	\boldsymbol{\varphi}_{k+1}(\theta_v,\theta_w) \   	
			\boldsymbol{z}_{k+1}(\theta_v,\theta_w)  ]' $ we have
			\begin{align*}
			\left[  \begin{array}{c} 	
			\boldsymbol{\varphi}_{k+1}(\theta_v,\theta_w) \\ 	
			\boldsymbol{z}_{k+1}(\theta_v,\theta_w) \end{array} \right] 
			& = B  \left[  \begin{array}{c} 	\boldsymbol{\varphi}_{k}(\theta_v,\theta_w) \\ 	\boldsymbol{z}_{k}(\theta_v,\theta_w) \end{array} \right] +  \left[  \begin{array}{c} 	\mathcal{L}_{k+1}^{\theta_v,\theta_w} \\ 	\mathcal{L}_{k+1}^{\theta_v,\theta_w} \end{array} \right] 
			\end{align*}
			which implies that for all $k\geq 1$,
			\begin{align*}
			& \left[  \begin{array}{c} 	\boldsymbol{\varphi}_{k+1}(\theta_v,\theta_w) \\ 	\boldsymbol{z}_{k+1}(\theta_v,\theta_w) \end{array} \right]  = B^{k+1} \left[  \begin{array}{c} 	\boldsymbol{\varphi}_{0}(\theta_v,\theta_w) \\ 	\boldsymbol{z}_{0}(\theta_v,\theta_w) \end{array} \right] \\
			&\qquad \qquad + \sum\limits_{t=1}^{k} B^{k+1-t}\left[  \begin{array}{c} 	\mathcal{L}_{t}^{\theta_v,\theta_w} \\ 	\mathcal{L}_{t}^{\theta_v,\theta_w}  \end{array} \right] + \left[  \begin{array}{c} 	\mathcal{L}_{k+1}^{\theta_v,\theta_w} \\ \mathcal{L}_{k+1}^{\theta_v,\theta_w} \end{array} \right] .
			\end{align*}
			
			Then we have
			\begin{align*}
			&\boldsymbol{\varphi}_{k}(\theta_v,\theta_w)  = [I_n \ \boldsymbol{0}]B^{k} [I_n  \ I_n ]'\boldsymbol{\varphi}_{0}(\theta_v,\theta_w) \\
			& \qquad+ \sum\limits_{t=1}^{k} [I_n \ \boldsymbol{0}]B^{k-t} [I_n \ I_n]'\mathcal{L}_{t}^{\theta_v,\theta_w}
			\end{align*}
			where the assumption of uniform initial beliefs sets the first term of the above relation to zero.
			
			The remainder of the proof follows the structure of the proof of Theorem~\ref{teo2}, where we invoke 
			Lemma~\ref{analog_lemma3} instead of Lemma \ref{lemma_angelia}. First, we will find a bound for the expected value of $\boldsymbol{\varphi}_{k}(\theta_v,\theta_w)$ and later we will show this is of bounded variations. In this case, we have
				\begin{align*}
				\mathbb{E}[\varphi^i_{k}(\theta_v,\theta_w)]
				& =  - \sum\limits_{t=1}^{k} \sum\limits_{j=1}^{n}[ [I_n \ \boldsymbol{0}]B^{k-t} [I_n \ I_n]']_{ij} H^j(\theta_v,\theta_w) .
				\end{align*}

			By adding and subtracting $\sum\limits_{t=1}^{k} \sum\limits_{j=1}^{n}\frac{1}{n} H^j(\theta_v,\theta_w)$  we obtain
			\begin{align*}
			& \mathbb{E}[ \varphi^i_{k}(\theta_v,\theta_w)]
			= -  \frac{k}{n} \sum\limits_{j=1}^{n}H^j(\theta_v,\theta_w) \\
			& +\sum_{t=1}^{k} \sum\limits_{j=1}^{n} \left(   \frac{1}{n}  -[ [I_n \ \boldsymbol{0}]B^{k-t} [I_n \ I_n]']_{ij}\right) H^j(\theta_v,\theta_w) .
			\end{align*}
			
			Similarly, as in the proof of Theorem \ref{teo2}, we bound the term in parenthesis using the non-asymptotic bounds from Lemma \ref{analog_lemma3} in conjunction with Lemma \ref{shahin}. By doing so, it can be seen that
			\begin{align*}
			\mathbb{E}[ \varphi^i_{k}(\theta_v,\theta_w)]
			& \leq  \frac{4 \log n}{1 - \lambda} \log\frac{1}{\alpha} - \frac{k}{n} \sum\limits_{j=1}^{n} H^j(\theta_v,\theta_w).
			\end{align*}
			
			Now, we will show that $ \varphi_{k}^i(\theta_v,\theta_w)$, as a function of the random variables consisting in $S_t^j$ for $1\leq t \leq k$ to $1\leq j\leq n$,  has bounded variations and we will compute the bound. First, we fix all other input random variables but $ [\mathcal{L}^{\theta_v,\theta_w}_t]_j$ and we have
			\begin{align*}
			& \max_{s_t^j \in \mathcal{S}^j} \varphi_{k}^i(\theta_v,\theta_w) - \min_{s_t^j \in \mathcal{S}^j} \varphi_{k}^i(\theta_v,\theta_w) 
			\\
			& =  \max_{s_t^j \in \mathcal{S}^j} 
			[ [I_n \ \boldsymbol{0}]B^{k-t} [I_n \ I_n]']_{ij}[ \mathcal{L}^{\theta_v,\theta_w}_t]_{j}  \\
			& \qquad -\min_{s_t^j \in \mathcal{S}^j}[ [I_n \ \boldsymbol{0}]B^{k-t} [I_n \ I_n]']_{ij}[\mathcal{L}^{\theta_v,\theta_w}_t]_{j}\\
			& \leq [ [I_n \ \boldsymbol{0}]B^{k-t} [I_n \ I_n]']_{ij}2\log\frac{1}{\alpha}.
			\end{align*}
			Thus, the summation of the squared bounds in McDiarmid's inequality is
			\begin{align*}
			&\sum_{t=1}^{k}\sum_{j=1}^{n} \left(  [ [I_n \ \boldsymbol{0}]B^{k-t} [I_n \ I_n]']_{ij} 2\log\frac{1}{\alpha}\right) ^2.
			\end{align*}
			
			Now, by adding and subtracting the term $1/n$ we have that
			\begin{align*}
			&\sum_{t=1}^{k}\sum_{j=1}^{n} \left(  [ [I_n \ \boldsymbol{0}]B^{k-t} [I_n \ I_n]']_{ij} 2\log\frac{1}{\alpha}\right) ^2 \\
			& \qquad \leq 8\left( \log\frac{1}{\alpha}\right) ^2\sum_{t=1}^{k}\sum_{j=1}^{n} \left(  [ [I_n \ \boldsymbol{0}]B^{k-t} [I_n \ I_n]']_{ij} - 1/n\right) ^2 \\
			& \qquad \qquad + 8\left( \log\frac{1}{\alpha}\right) \sum_{t=1}^{k}\sum_{j=1}^{n} \left( 1/n\right) ^2 
			\end{align*}
			where we have used $x^2 \leq 2((x-y)^2+y^2)$.
			\newline
			We can bound the first term in the preceding relation using Eq. \eqref{square_bound} with $\boldsymbol{y}_1 = \boldsymbol{e}_j$ since Eq. \eqref{square_bound} hold for any choice of $\boldsymbol{y}_1$. Specifically, we obtain that for all $j = 1,\hdots ,n$
			\begin{align*}
\sum_{i=1}^n \left(  [[I_n \ \boldsymbol{0}]B^{k-t} [I_n  \ I_n ]']_{ij} - 1/n \right)^2\leq 2 \left(1 - \frac{1}{9U}\right)^{k-t}.
\end{align*}
\newline
Additionally, note that $[I_n \ \boldsymbol{0}]B^{k} [I_n  \ I_n ]'$ is a symmetric matrix since is a polynomial of $\bar{A}$ which is symmetric itself. This in turn implies that $[I_n \ \boldsymbol{0}]B^{k-t} [I_n  \ I_n ]'$ is also symmetric.  Therefore, it holds that for all $i = 1,\hdots ,n$
			\begin{align*}
\sum_{j=1}^n \left(  [[I_n \ \boldsymbol{0}]B^{k-t} [I_n  \ I_n ]']_{ij} - 1/n \right)^2\leq 2 \left(1 - \frac{1}{9U}\right)^{k-t}.
\end{align*}

Finally, we have
			\begin{align*}
&\sum_{t=1}^{k}\sum_{j=1}^{n} \left(  [ [I_n \ \boldsymbol{0}]B^{k-t} [I_n \ I_n]']_{ij} 2\log\frac{1}{\alpha}\right) ^2 \\
& \qquad \leq 8\left( \log\frac{1}{\alpha}\right) ^2\left( 2k + \frac{k}{n}\right)  \\
			& \qquad \leq 24\left( \log \alpha\right) ^2 k.
			\end{align*}
						
			Now, by the McDiarmid inequality and getting the values of $k$ such that the desired probabilistic tolerance level $\rho$ is achieved, we obtain
			\begin{align*}
			&\mathbb{P}\left(\varphi_{k}^i(\theta_v,\theta_w) - \mathbb{E}[\varphi_{k}^i(\theta_v,\theta_w) ] \geq   \frac{k}{2}\gamma_2  \right) \\
			& \qquad  \qquad=  \exp\left( -\frac{2\left(\frac{1}{2} k\gamma_2 \right)^2}{ 24\left( \log \alpha\right) ^2 k}  \right) \\
			& \qquad \qquad =  \exp\left( -\frac{k\gamma_2 ^2}{ 48\left( \log \alpha\right) ^2}  \right). 
			\end{align*}
			
			Therefore, for a given confidence level $\rho$, in order to have 
			$\mathbb{P}\left(\mu_{k}^i(\theta_v) \geq \exp\left( -\frac{1}{2}k\gamma_2 + \gamma_1^i \right) \right) 
			\leq \rho$ we require that
			\begin{align*}
			k & \geq \frac{1}{\gamma_2^2} 48 \left( \log\alpha \right) ^2 \log\frac{1}{\rho} .
			\end{align*}
		\end{IEEEproof}

\begin{figure}[ht]
			\centering
			\subfigure[Path Graph]{
				\includegraphics[width=0.22\textwidth]{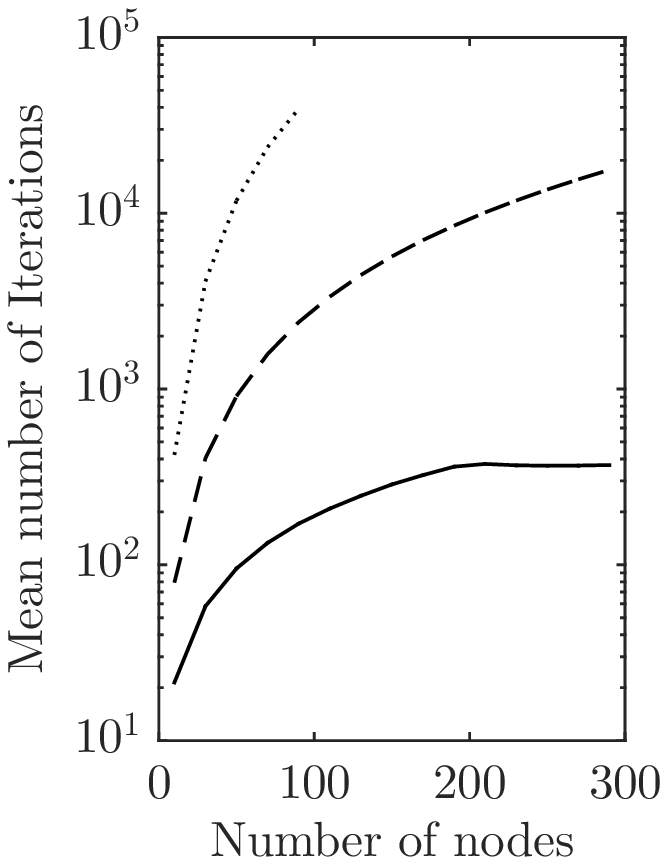}
			}
			\subfigure[Circle Graph]{
				\includegraphics[width=0.22\textwidth]{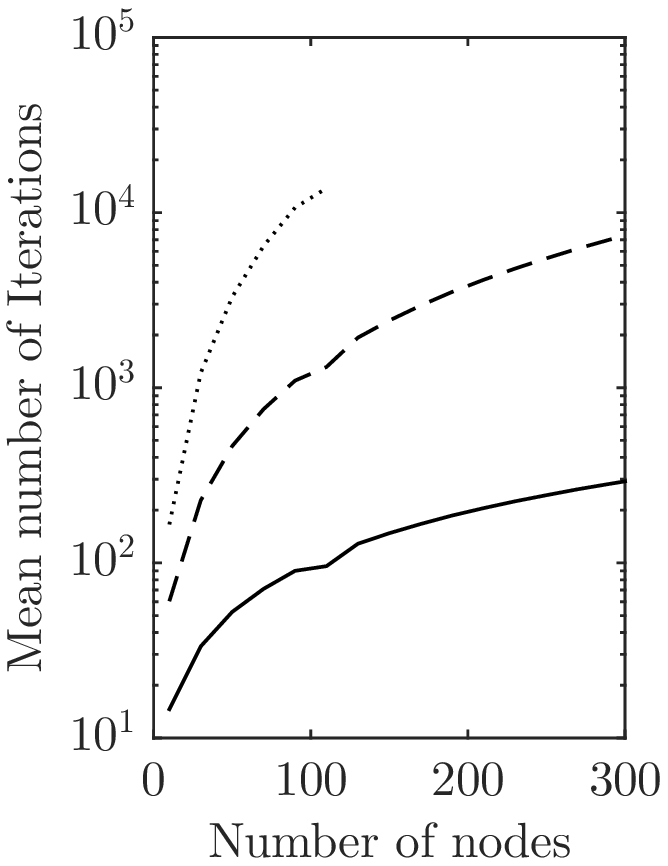}
			}
			\\
			\subfigure[Grid Graph]{
				\includegraphics[width=0.22\textwidth]{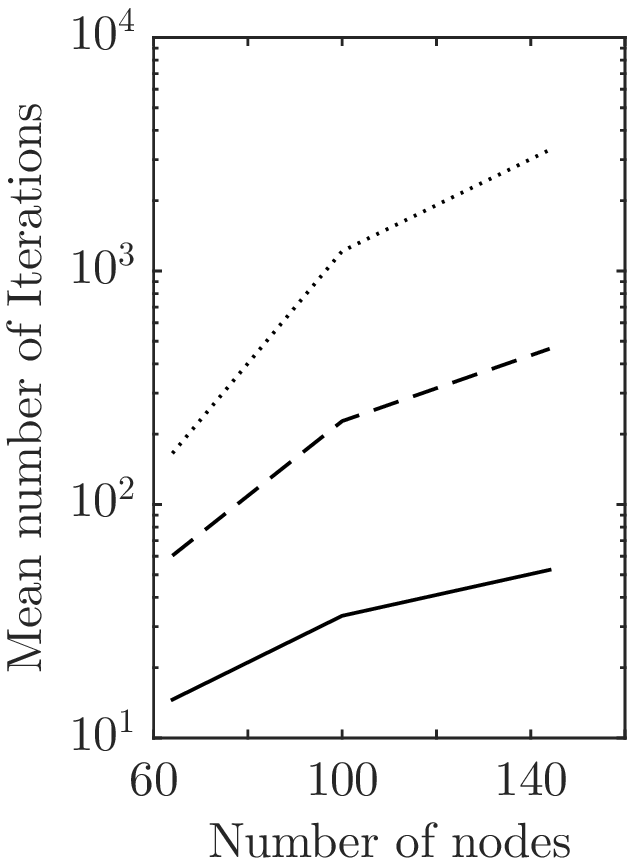}
			}
			\caption{Empirical mean over 50 Monte Carlo runs of the number of iterations required for $\mu_k^i(\theta) < \epsilon$ for all agents on $\theta \notin \Theta^*$. All agents but one have all their hypotheses to be observationally equivalent. Dotted line for the algorithm proposed in \cite{jad12}, Dashed line for the procedure described in Eq.~\eqref{non_bayes_distributed} and solid line for the procedure described in Eq.~\eqref{linear_bayes}.}
			\label{one_learn}
		\end{figure}

		Next, we will present simulation results that show how the convergence time depends on the number of agents in the network. Figure~\ref{one_learn} shows the time required for a group of agents to have a set of beliefs at a distance of $\epsilon = 0.01$ from the singleton distribution around the optimal hypothesis. For example, on a path graph, as the path grows longer, the number of iterations required to meet the desired $\epsilon$ accuracy grows rapidly. This is due to the low connectivity of the network. The time required for consensus is smaller for the circle and the grid graphs due to their better connectivity properties.
		
		\section{Numerical Example: Distributed Source Localization}\label{special}
        
In this section we apply the proposed algorithms to the problem of distributed source localization based on differential signal amplitudes \cite{rab04,rab05,wan12,mao07,lan03}. We compare the performance of our methods, Eq. \eqref{non_bayes_distributed} and Eq. \eqref{linear_bayes} with the algorithms proposed in \cite{qip15,jad12}. For simulation purposes we will assume the graphs are fixed and there exists a single $\theta^*$ such that $f^i = \ell^i\left(\cdot|\theta^*\right)$ for all $i$, in which case our update rule simplifies to the learning algorithm proposed in \cite{sha14}.

\begin{figure}[ht]
			\centering
			\subfigure[Network of $3$ agents]{
				\includegraphics[width=0.35\textwidth]{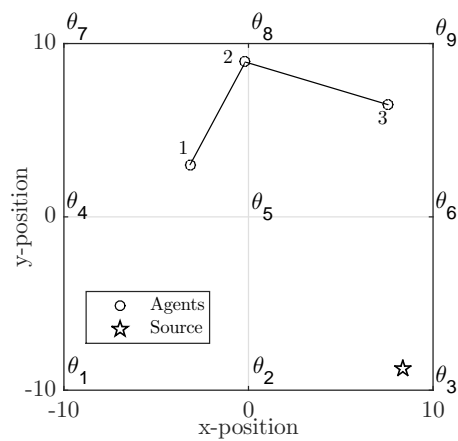}
			}
			\subfigure[Distributions for Hypothesis $\theta_5$ and $\theta_3$.]{
				\includegraphics[width=0.4\textwidth]{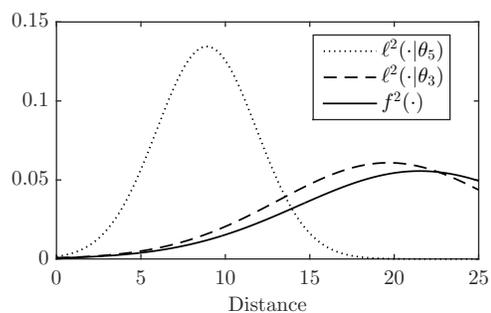}
			}
			\caption{(a) Group of 3 agents in a grid of $3 \times 3$ hypotheses. Each hypothesis corresponds to a possible location of the source. For example, hypothesis $\theta_2$ locates the source at the $\left(-10,0\right)$ point in the plane. (b) Likelihood functions for $\theta_2$ and $\theta_5$ and distribution of observations $f^2$ for agent 2.}
			\label{grid_target}
		\end{figure}

Assume a group of $n$ agents is randomly distributed in an area and each agent receives a noisy signal proportional to its distance to a target. The group objective is to collectively find the location of the target. Each agent constructs a grid of hypotheses about the possible location of the source. Figure \ref{grid_target}(a) shows a $10$ by $10$ area partitioned in a $3$ by $3$ grid, which results in $9$ hypotheses. Moreover, there are three agents (represented by circles), at different locations. The graph structure shows that agent $1$ communicates with agent $2$, similarly, agent $3$ communicates with $2$. The star represents the target.

Each agent constructs likelihood functions for its hypotheses based on its sensor model. The observations follow a truncated normal distribution with the mean proportional to the distance between the agent and the grid point of the corresponding hypothesis. For example, assume an agent $i$ is in a position $p_a^i = (x_a^i,y_a^i)$ and the target is located at $p_s = (x_s,y_s)$. The received signals are $S_k^i = \|p_s - p_a^i\| + c W_k^i$, where $c$ is some positive constant and $W_k^i$ is a truncated zero mean Gaussian noise. Now, consider that a hypothesis $\theta$ is at a point $p_\theta = (x_\theta,y_\theta)$. The corresponding likelihood model under hypothesis $\theta$ assumes observations are ${S_k^i|\theta = \|p_\theta - p_a\| + cW_k^i}$.

Figure \ref{grid_target}(b) shows the likelihood functions for $\theta_5$ and $\theta_3$ of agent 2, clearly hypothesis $\theta_3$ is closer to the true distribution of the observations $f^2$. Note that there is not a ``true state of the world" in the sense that $f^2$ is not equal to any of the hypotheses in the grid.

\begin{figure}[ht]
			\centering
			\includegraphics[width=0.5\textwidth]{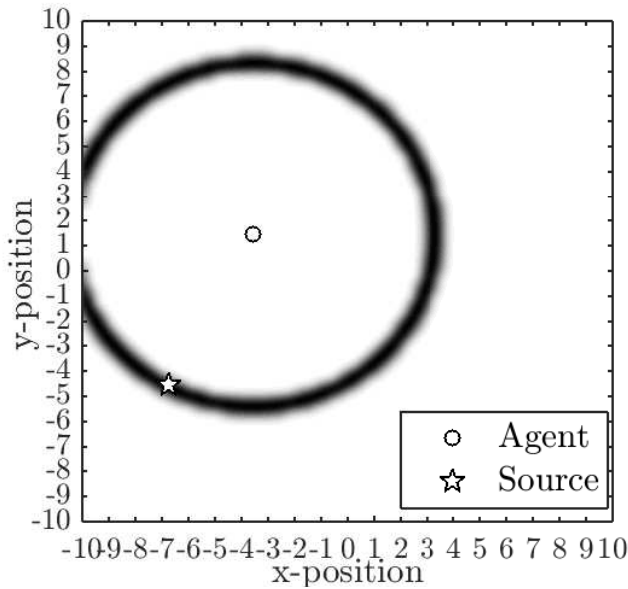}
			\caption{Belief distribution of one agent over the hypotheses grid. Darker shades of gray indicates higher beliefs on the corresponding hypothesis.}
			\label{single_agent_circle}
		\end{figure}

		The information each agent obtains is enough just to estimate the distance to the source, but not its complete coordinates. For instance, a single sensor can only locate the source within a circular band around it, see Figure \ref{single_agent_circle}. 
		
\begin{figure}[ht]
			\centering
			\subfigure[Network of Agents]{
				\includegraphics[width=0.4\textwidth]{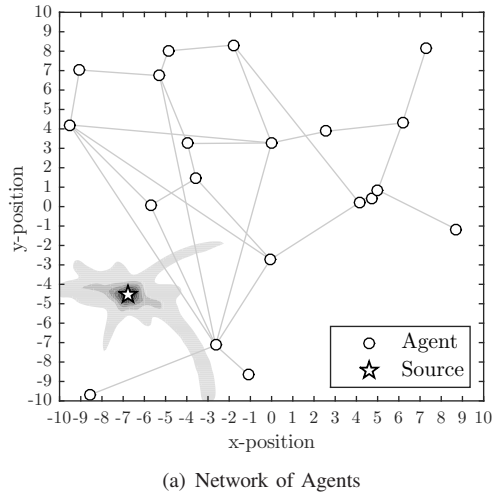}
			}
			\subfigure[Belief of one agent on the optimal hypothesis]{
				\begin{overpic}[width=0.4\textwidth]{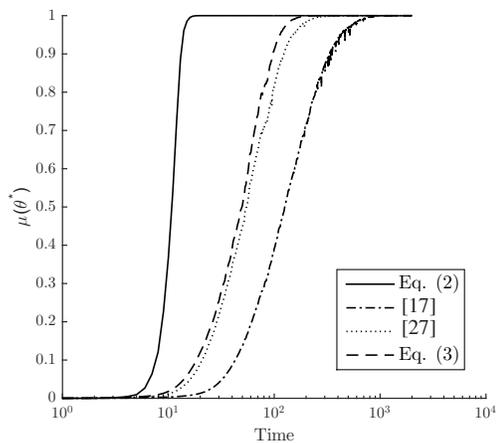}
					\put(74,19){\color{white}{\rule{0.5cm}{0.7cm}}}
					\put(75,25){{{\footnotesize  \cite{jad12}}}}
					\put(75,21){{{\footnotesize  \cite{qip15}}}}
				\end{overpic} 
			}
			\caption{(a) Network of agents as well as the belief distribution over the hypothesis set (a grid in the $x$, $y$ location). Darker shade of gray indicates higher beliefs on the corresponding hypothesis (point in the hypotheses grid). (b) Belief evolution on the optimal hypothesis $\theta^*$ for different belief update protocols.}
			\label{learning_network}
		\end{figure}

		Figure \ref{learning_network}(a) shows another group of $20$ agents now interacting according to an appropriate network structure, see Assumptions \ref{assum:graph} and \ref{assum_linear}. A finer grid partition has been used, where each coordinate has $100$ points, resulting in $10000$ hypotheses in total. Figure \ref{learning_network}(b) shows the belief on the hypothesis $\theta^*$, defined to be the grid point closer to the location of the target.
		
\begin{figure}[ht]
			\centering
			\subfigure[Network of agents]{
				\includegraphics[width=0.4\textwidth]{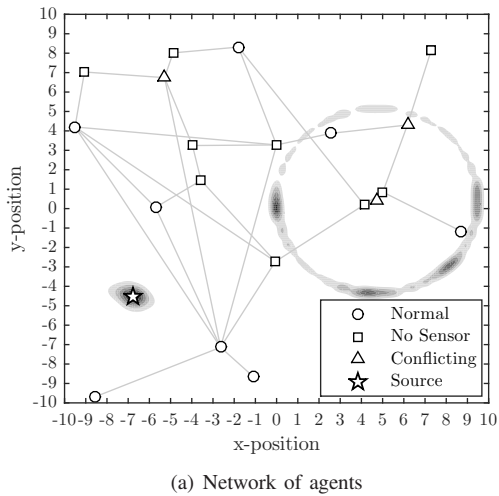}
			}
			\subfigure[Belief of one agent on the optimal hypothesis]{
				\begin{overpic}[width=0.4\textwidth]{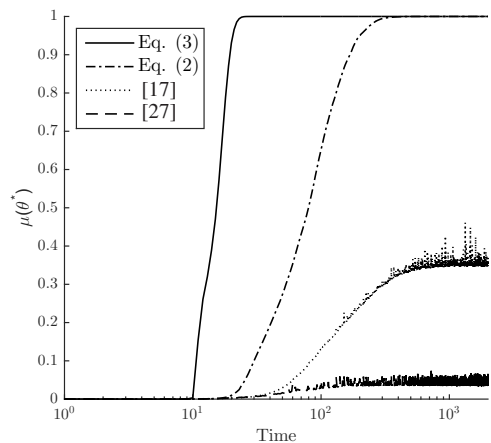}
					\put(26,60){\color{white}{\rule{0.5cm}{0.6cm}}}
					\put(27,65){{{\footnotesize  \cite{jad12}}}}
					\put(27,61){{{\footnotesize  \cite{qip15}}}}
				\end{overpic} 
			}
			\caption{(a) Network of heterogeneous agents. $\triangle$ indicates agents whose observations have been modified such that the optimal hypothesis is the $\left(0,0\right)$ point in the grid. $\square$ indicates agents for whom all hypotheses are observationally equivalent (i.e. no data is measured). $\circ$ indicates regular agents with correct observation models and informative hypothesis. (b) Belief evolution on the optimal hypothesis $\theta^*$ for different belief update protocols. }
			\label{learning_network_conflict}
		\end{figure}

		Figure \ref{learning_network_conflict} repeats the simulations presented in Figure \ref{learning_network} but including $10$ agents with all their hypotheses observationally equivalent (i.e. no measurements available), and $3$ conflicting agents whose observations have been modified (corrupted) such that the optimal hypothesis is the $\left(0,0\right)$ point in the grid. 
		
		Figure \ref{learning_network_conflict}(b) shows the protocols presented in Eqs. \eqref{non_bayes_distributed} and \eqref{linear_bayes} concentrate the beliefs onto the optimal hypothesis. The performance of the algorithms in \cite{jad12} and \cite{qip15} deteriorates if conflicting agents are present. This is evident from the lack of concentration of the beliefs around the true hypotheses.
		
		\section{Conclusions and Future Work}\label{conclusions}
		We proposed two distributed cooperative learning algorithms for the problem of collaborative inference. The first algorithm focuses on general time-varying undirected graphs, and the second algorithm is specialized for fixed graphs. In both cases, we show that the beliefs converge to the hypothesis set that best describes the observations in the network. We require reasonable connectivity assumptions on the communication network over which the agents exchange information. 
		
		Our results prove convergence rates that are non-asymptotic, geometric, and explicit. The bounds depend explicitly on the graph sequence properties, as well as the agent learning capabilities. Moreover, we do so in a new general setting where there might not be a ``true state of the world" which is perfectly described by a single hypothesis, i.e. misspecified models. Additionally, we analyze networks where agents might have conflicting hypotheses, i.e. the hypotheses with the highest confidence changes if different subsets of agents are taken into account. The algorithm for fixed undirected graphs achieves a factor of $n$ improvement in the convergence rate with respect to the number of agents in comparison with that of the existing algorithms.
		
		Our work suggests a number of open questions. The problem of tracking optimal hypothesis when its distributions are changing with time requires further study \cite{bat14}. Ideas from social sampling can also be incorporated in this framework \cite{sar13}, where the dimension of the beliefs is large and only partial beliefs are transmitted. Moreover, studying the influence of corrupted measurements or malicious agents is also of interest, especially in the setting of social networks. 
		
		%
		
		\bibliographystyle{IEEEtran} 
		
		\bibliography{IEEEfull,bayes_cons_3}
		
		\begin{IEEEbiography}[{\includegraphics[width=1in,height=1.25in,clip,keepaspectratio]{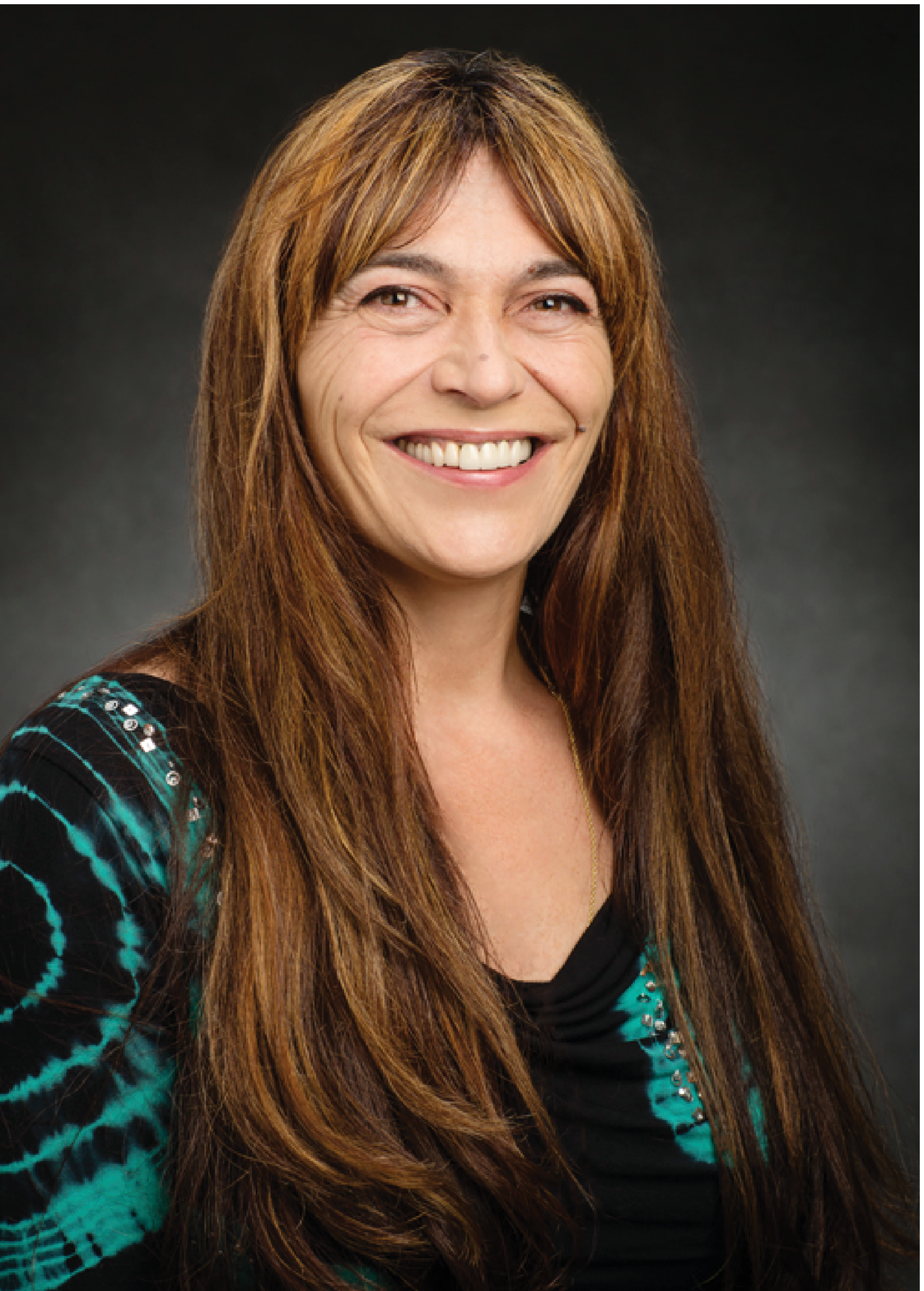}}]%
			{Angelia Nedi\'{c}}
			received the Ph.D. degree from Moscow State University, Moscow, Russia, in Computational Mathematics and Mathematical Physics, in 1994, and from Massachusetts Institute of Technology, Cambridge, USA, in electrical engineering and computer science, in 2002. 
			
			She is currently a Professor at the School of Electrical, Computer, and Energy Engineering, Arizona State University,Tempe, AZ, USA. Prior to this position, she served as a Faculty Member at the University of Illinois in Urbana-Champaign, from 2006 to 2016. She has worked as a Senior Engineer in BAE Systems North America, Advanced Information Technology Division at Burlington, MA. Her research interest is in large-scale optimization, games, control and information processing in networks.
			
			Prof. Nedi\'{c} received the NSF CAREER Award 2007 in operations research for her work in distributed multiagent optimization, the Best Paper Award at the Winter Simulation Conference 2013 (jointly with F. Yousefian and U. V. Shanbhag), and the Best Paper Award at the International Symposium on Modeling and Optimization in Mobile, Ad Hoc and Wireless Networks 2015 (jointly with K. Cohen and R. Srikant). 
			
		\end{IEEEbiography}
		\begin{IEEEbiography}[{\includegraphics[width=1in,height=1.25in,clip,keepaspectratio]{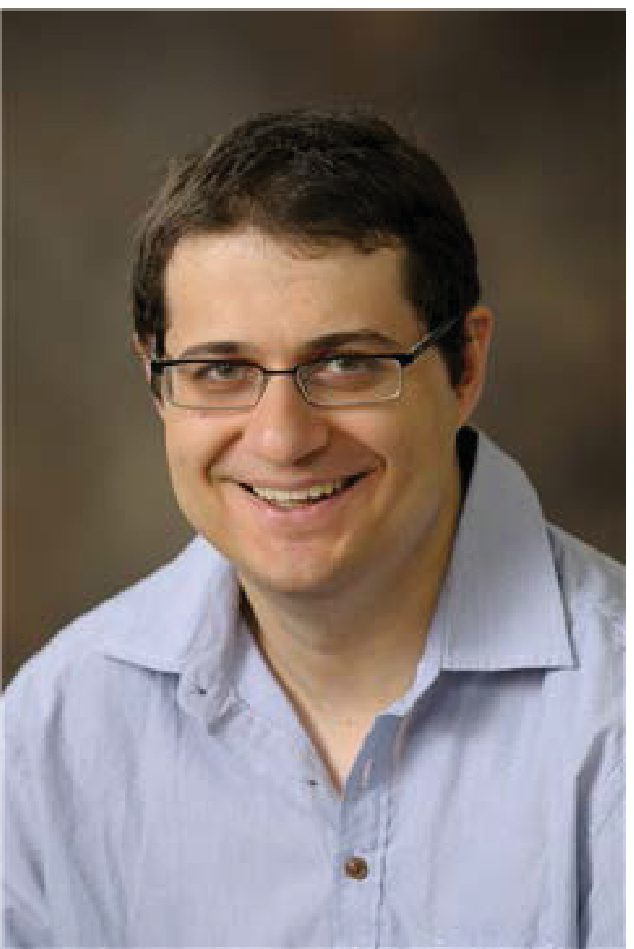}}]%
			{Alex Olshevsky}
			received the B.S. degree in applied
			mathematics and in electrical engineering
			from the Georgia Institute of Technology, Atlanta, GA, USA, both in 2004, and the M.S. and Ph.D. degrees in
			electrical engineering and computer science from the
			Massachusetts Institute of Technology, Cambridge, MA, USA, in 2006 and
			2010, respectively. 
			
			He is currently an Assistant Professor at the Department of Electrical and Computer Engineering, Boston University, Boston, MA, USA. Prior to this position, he served as a Faculty Member at the University of Illinois at Urbana-Champaign. He was a Postdoctoral Scholar in the Department of Mechanical and Aerospace Engineering, Princeton
			University, from 2010 to 2012 before joining the University of Illinois at
			Urbana-Champaign in 2012. His research interests are in control theory and
			optimization, especially in distributed, networked, and multi-agent settings.
			
			Dr. Olshevsky received the NSF CAREER Award, the Air Force Young
			Investigator Award, the ICS Prize from INFORMS for best paper on the
			interface of operations research and computer science, and a SIAM Paper
			Prize for annual paper from the \textit{SIAM Journal on Control and Optimization}
			chosen to be reprinted in SIAM Review.
		\end{IEEEbiography}
		\begin{IEEEbiography}[{\includegraphics[width=1in,height=1.25in,clip,keepaspectratio]{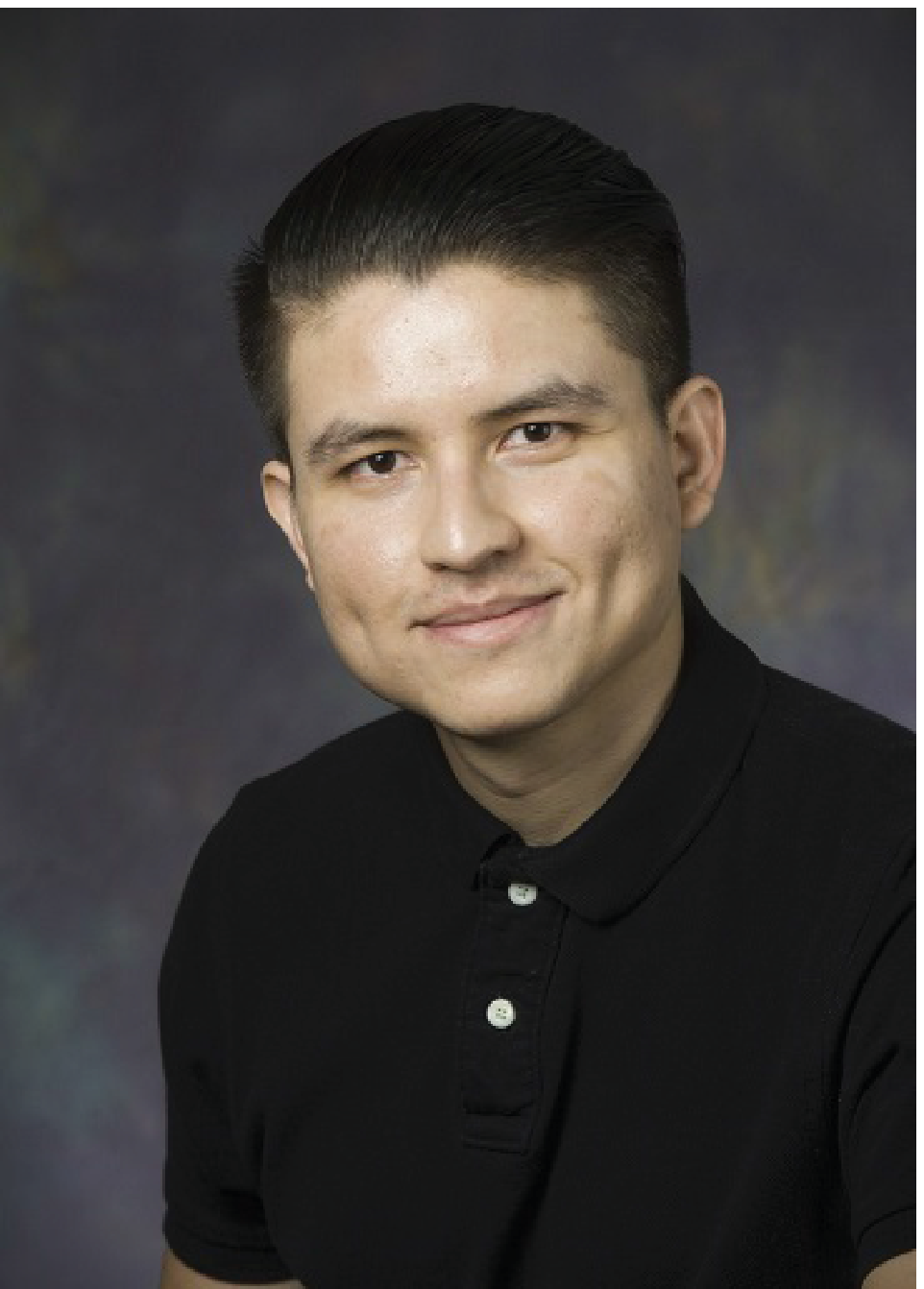}}]%
			{C\'{e}sar A. Uribe}
			received the B.S. degree in electronic engineering from the University of Antioquia, Medell\'{i}n, Colombia, in 2010, the M.S. degrees in systems and control from Delft University of Technology, Delft, Netherlands and in applied mathematics from the University of Illinois at Urbana-Champaign, Champaign, IL, USA, in 2013 and 2016, respectively. He is currently working toward the Ph.D. degree in electrical and computer engineering at the University of Illinois at Urbana-Champaign.
			
			His research interests include distributed learning and optimization, decentralized control, optimization theory, and algorithm analysis.
		\end{IEEEbiography}
	\end{document}